\documentclass[final,onefignum,onetabnum]{article}

\usepackage{array}
\usepackage{enumitem}
\usepackage{esvect}

\usepackage[caption=false]{subfig}
\usepackage{pgfplots}

\usepackage[colorlinks=true, bookmarksopen,
            pdfsubject={algorithms},
            linkcolor={blue},
            anchorcolor={black},
            citecolor={black},
            filecolor={magenta},
            menucolor={black},
            plainpages=false,pdfpagelabels,
            urlcolor={blue}]{hyperref}

\usepackage{algorithmic}

\usepackage{graphicx,epstopdf}

\usepackage{subfig}
\usepackage{amsmath}
\usepackage{amssymb}
\usepackage{graphicx}
\usepackage{dcolumn}
\usepackage{mathtools}
\usepackage{amsmath,amsfonts}
\usepackage{bm}
\usepackage{amsmath}
\usepackage{amssymb}
\usepackage{color}
\usepackage{float}
\usepackage{setspace}
\usepackage{tabularx}

\allowdisplaybreaks

\newcommand{\f}{\frac}

\newcommand{\E}{ {\mathbb{E}} }
\newcommand{\ola}{\overleftarrow}

\newcommand{\be}{\begin{equation}}
\newcommand{\ee}{\end{equation}}

\def\ba{\begin{array}}                \def\ea{\end{array}}
\def\bel{\begin{equation}\label}      \def\ee{\end{equation}}

\colorlet{texcscolor}{blue!50!black}
\colorlet{texemcolor}{red!70!black}
\colorlet{texpreamble}{red!70!black}
\colorlet{codebackground}{black!25!white!25}

\date{}

\begin{document}


\title{A Score-based Nonlinear Filter for Data Assimilation }

%
%

%
%
%
\author{
Feng Bao\thanks{ Department of Mathematics, Florida State University, Tallahassee, Florida, \ ({\tt bao@math.fsu.edu}).} 
\and Zezhong Zhang \thanks{ Department of Mathematics, Florida State University, Tallahassee, Florida.}
\and Guannan Zhang  \thanks{ Computer Science and Mathematics Division, Oak Ridge National Laboratory, Oak Ridge, TN 37831, USA. 
 This manuscript has been authored by UT-Battelle, LLC, under contract DE-AC05-00OR22725 with the US Department of Energy (DOE). The US government retains and the publisher, by accepting the article for publication, acknowledges that the US government retains a nonexclusive, paid-up, irrevocable, worldwide license to publish or reproduce the published form of this manuscript, or allow others to do so, for US government purposes. DOE will provide public access to these results of federally sponsored research in accordance with the DOE Public Access Plan.}
      }
      
      
 \maketitle
      
\begin{abstract}
We introduce a score-based generative sampling method for solving the nonlinear filtering problem with robust accuracy. A major drawback of existing nonlinear filtering methods, e.g., particle filters, is the low stability. To overcome this issue, we adopt the diffusion model framework to solve the nonlinear filtering problem. In stead of storing the information of the filtering density in finite number of Monte Carlo samples, in the score-based filter we store the information of the filtering density in the score model. Then, via the reverse-time diffusion sampler, we can generate unlimited samples to characterize the filtering density. Moreover, with the powerful expressive capabilities of deep neural networks, it has been demonstrated that a well trained score in diffusion model can produce samples from complex target distributions in very high dimensional spaces. Extensive numerical experiments show that our score-based filter could
potentially address the curse of dimensionality in very high dimensional problems.  
\end{abstract}

\textbf{keyword:}
Nonlinear filtering, diffusion model, score-based models, stochastic dynamical systems.

\section{Introduction}

Nonlinear filtering is a major research direction in data assimilation, which has a broad spectrum of applications in weather forecasting, military, material sciences, biology, and finance \cite{Bao_CiCP20, Bao_Cogan20, Multi-target_2007, Bao_Atom20, Evense_EnKF, Filter_Finance}. The goal of solving a filtering problem is to exploit partial noisy observational data stream to estimate the unobservable state of a stochastic dynamical system of interest. 
%
%
%
%
In linear filtering, i.e., both the state dynamics and the observation dynamics are linear, the optimal estimate for the unobservable state can be analytically obtained by the Kalman filter under the Gaussian assumption. 

When the dynamical systems are nonlinear, the standard Kalman filter is no longer feasible. An extension of the Kalman filter, i.e., the ensemble Kalman filter, can be applied to address the nonlinearity to some extent. The main idea of the ensemble Kalman filter is to use an ensemble of Kalman filter samples to characterize the probability distribution of the target state as a Gaussian distribution. As a Kalman type filter, the ensemble Kalman filter stores the information of the state variable as the mean and the covariance of Kalman filter samples. In this way, the probability density function (PDF) of the target state, which is often called the ``filtering density'', is approximated as a Gaussian distribution.  However, in nonlinear filtering problems, the filtering density is usually non-Gaussian. Therefore, the ensemble Kalman filter, which still relies on the Gaussian assumption, is not the ideal approach to solve the nonlinear filtering problem.

Besides the ensemble Kalman filter, several effective nonlinear filtering methods (e.g., the particle filter \cite{MCMC-PF, particle-filter}, the Zakai filter \cite{Bao_Zakaid_2015, zakai}, and the backward SDE filter \cite{Bao_AA20, BaoC20142, Bao_first, BSDE_filter})  are developed to overcome the nonlinearity. Among those methods, the particle filter is the most widely applied approach for solving the nonlinear filtering problem. The particle filter method is also known as the ``sequential Monte Carlo'' method. It utilizes a set of Monte Carlo samples, called ``particles'', to construct an empirical distribution to describe the filtering density of the target state. Upon reception of the observational data, a Bayesian inference procedure is adopted to assign likelihood as weights to those particles, and a resampling procedure is repeatedly carried out to generate more duplications of higher weighted particles and drop particles with low weights.  The particle filter can handle the nonlinearity issue in the optimal filtering problem since the nonlinear state dynamics can be incorporated into the filtering density through particle simulations, and the Bayesian inference is the standard approach to deal with nonlinear observations.  Different from Kalman type filters, which store the information of the filtering density only in the mean and the covariance, the particle filter stores the information of the filtering density in locations of various particles. Therefore, the particle filter can characterize more complex non-Gaussian filtering densities via empirical distributions built by the particles.

The main drawback of the particle filter is its low stability. In the particle filter, once the total number of particles is fixed, the capability to characterize the filtering density is fixed, and one may only use those finite particles to approximate the filtering density.  However, the nonlinear state dynamics and nonlinear observations could result unpredictable features for filtering densities. Therefore, it's hard to use finite particles to characterize unlimited possibilities of filtering densities. This often causes the so-called ``degeneracy issue'', i.e. the particles lie in high probability regions are not sufficient to characterize highly probable features in filtering densities. Such a degeneracy issue is even more prohibitive when the dimension of the problem is high due to ``curse of dimensionality''. Although advanced resampling methods are proposed to address the degeneracy issue by relocating particles to high probability regions \cite{Kang-PF, MCMC-PF, APF, Sny, CT1}, the nature of finite particle representation for filtering density cannot be changed under the sequential Monte Carlo framework, and the information that any resampling method may use cannot exceed the information carried by those finite particles.


In this work, we introduce a novel score-based filter method that allows to use a ``score'' model to generate samples from the filtering density through a diffusion process. The \textit{score} is a key component in diffusion model, which is a well-known generative machine learning model for generating samples from a target PDF. 
Diffusion models are generative models that utilize noise injection to progressively distort data and then learn to reverse this process for sample generation. 
As a category of deep generative models, diffusion models are widely used in image processing applications, such as image synthesis \cite{NEURIPS2021_49ad23d1,NEURIPS2020_4c5bcfec,NEURIPS2019_3001ef25,song2021scorebased,DBLP:conf/eccv/CaiYAHBSH20,DBLP:journals/jmlr/HoSCFNS22,DBLP:conf/iclr/MengHSSWZE22}, image denoising \cite{DBLP:conf/iccvw/KawarVE21,DBLP:conf/iccv/LuoH21,NEURIPS2020_4c5bcfec,DBLP:conf/icml/Sohl-DicksteinW15}, image enhancement \cite{DBLP:journals/corr/abs-2112-05149,DBLP:journals/ijon/LiYCCFXLC22,DBLP:journals/pami/SahariaHCSFN23,DBLP:conf/cvpr/WhangDTSDM22}, image segmentation \cite{amit2022segdiff,baranchuk2022labelefficient,DBLP:conf/cvpr/BrempongKCPM022,DBLP:journals/corr/abs-2206-09012}, and natural language processing \cite{DBLP:conf/nips/AustinJHTB21,DBLP:conf/nips/HoogeboomNJFW21,DBLP:journals/corr/abs-2205-14217,DBLP:conf/iclr/SavinovCBEO22,DBLP:conf/icml/YuXMJPGZZW22}. 
Diffusion models are also capable of density estimation (i.e., learning how to draw samples from the probability distribution). Specifically,
a diffusion model can transport a prior distribution, which is often chosen as the standard Gaussian distribution, to a complex target data distribution through a reverse-time diffusion process in the form of a stochastic differential equation, and the score model is the ``forcing term'' that guides the diffusion process towards the data distribution. Since the prior distribution is independent from the target data distribution, the information of the data distribution is stored in the score model. When adopting the diffusion model framework to solve the nonlinear filtering problem, we let the filtering density be our target distribution. Specifically, we propagate Monte Carlo samples through the state dynamics to generate ``data samples'' that follow the filtering density, and we use those data samples to train a score model in the form of a deep neural network. Although samples that characterize the filtering density are still needed in our method, the score-based filter is essentially different from current Monte Carlo based filtering methods. In fact, in stead of storing the information of the filtering density in finite number of Monte Carlo samples, in the score-based filter we store the information of the filtering density in the score model. Then, via the reverse-time diffusion sampler, we can generate unlimited samples to characterize the filtering density. Moreover, with the powerful expressive capabilities of deep neural networks, it has been demonstrated that a well trained score in diffusion model can produce samples from complex target distributions in very high dimensional spaces \cite{song2021scorebased}. Therefore, the score-based filter that we develop in this paper may provide a promising numerical method that could potentially address the curse of dimensionality in very high dimensional problems. 

The rest of this paper is organized as follows. In Section \ref{Problem setting}, we briefly introduce the nonlinear filtering problem and its state-of-the-art solver, i.e. the particle filter method. In Section \ref{sec:SNF}, we provide a comprehensive discussion to develop our score-based filter method. In Section \ref{sec:exam}, we carry out numerical experiments for three benchmark nonlinear filtering problems, which include a $100$-dimensional Lorenz-96 attractor problem, to show the high robustness and high accuracy of the score-based filter.

\section{Problem setting}\label{Problem setting}
Nonlinear filters are important tools for dynamical data assimilation with a variety of scientific and engineering applications. The definition of a nonlinear filtering problem can be viewed as an extension of Bayesian inference to the estimation and prediction of a nonlinear stochastic dynamical system. In this effort, we consider the following state-space nonlinear filtering model:
\begin{equation}\label{NLF}
\hspace{-2cm}
\begin{aligned}
\text{\small \bf State:}\quad & X_{t+1} =\; f(X_t, \omega_t),   \\
 \text{\small \bf  Observation:}\quad & Y_{t+1} \;=\;  g(X_{t+1})+ \varepsilon_{t+1}, 
\end{aligned}
\end{equation}
where $t \in \mathbb{Z}^{+}$ represents the discrete time, $X_t \in \mathbb{R}^d$ is a $d$-dimensional unobservable dynamical state governed by the nonlinear function $f: \mathbb{R}^d \times \mathbb{R}^k \mapsto \mathbb{R}^d$, $\omega_t \in \mathbb{R}^k$ is a random variable that follows a given probability law representing the uncertainty in $f$, and the random variable $Y_{t+1} \in \mathbb{R}^r$ provides nonlinear partial observation on $X_{t+1}$, i.e., $g(X_{t+1})$, perturbed by a Gaussian noise $\varepsilon_{t+1} \sim  \mathcal{N}(0, \Sigma)$. 

The overarching goal is to find the best estimate, denoted by $\hat{X}_{t+1}$, of the unobservable state $X_{t+1}$, given the observation data $\mathcal{Y}_{t+1}:= \sigma(Y_{1:t+1})$ that is the $\sigma$-algebra generated by all the observations up to the time instant $t+1$.
Mathematically, such optimal estimate for $X_{t+1}$ is usually defined by a conditional expectation, i.e.,
\begin{equation}\label{Expectation:S}
\hat{X}_{t+1} := \E[X_{t+1} | \mathcal{Y}_{t+1}],
\end{equation}
where the expectation is taken with respect to the random variables $\omega_{0:t}$ and $\varepsilon_{1:t+1}$ in Eq.~\eqref{NLF}. In practice, the expectation in Eq.~\eqref{Expectation:S} is  not approximated directly. Instead, we aim at approximating the conditional probability density function (PDF) of the state, denoted by
$P(X_{t+1} | \mathcal{Y}_{t+1})$,
which is referred to as the filtering density.
The Bayesian filter framework is to recursively incorporate observation data to describe the evolution of the filtering density. There are two steps from time $t$ to $t+1$, i.e., the prediction step and the update step:
\begin{itemize}[leftmargin=15pt]
\item {\em The prediction step} is to use the Chapman-Kolmogorov formula to propagate the state equation in Eq.~\eqref{NLF} from $t$ to $t+1$ 
and obtain the prior filtering density, i.e.,
\begin{equation}\label{Kolmogorov}
\hspace{-1cm}\text{\small \bf Prior filtering density:}\quad P(X_{t+1}\big|\mathcal{Y}_t) = \int P(X_{t+1} \big| X_t)P(X_{t}\big|\mathcal{Y}_t) dX_t,
\end{equation}
where $P(X_{t}\big|\mathcal{Y}_t)$ is the posterior filtering density obtained at the time instant $t$, $P(X_{t+1} \big| X_t)$ is the transition probability derived from the state dynamics in Eq.~\eqref{NLF}, and $P(X_{t+1}\big|\mathcal{Y}_t)$ is the prior filtering density for the time instant $t+1$.

\item {\em The update step} is to combine the likelihood function, defined by the new observation data $Y_{t+1}$, with the prior filtering density to obtain the posterior filtering density, i.e.,
\begin{equation}\label{Bayes}
\hspace{-1.5cm}\text{\small \bf Posterior filtering density:}\quad P(X_{t+1}\big|\mathcal{Y}_{t+1}) \propto {P(X_{t+1}\big|\mathcal{Y}_t) \, P(Y_{t+1} \big| X_{t+1})},
\end{equation} 
where the likelihood function $P(Y_{t+1} \big| X_{t+1})$ is defined by
\begin{equation}\label{Likelihood}
P(Y_{t+1} \big| X_{t+1}) \propto \exp\left[ -\f{1}{2} \big(g(X_{t+1}) - Y_{n+1}\big)^{\top}\Sigma^{-1} \big(g(X_{t+1}) - Y_{t+1} \big) \right],
\end{equation}
with $\Sigma$ being the covariance matrix of the random noise $\varepsilon$  in Eq.~\eqref{NLF}.

\end{itemize}

In this way, the filtering density is predicted and updated through formulas Eq.~\eqref{Kolmogorov} to Eq.~\eqref{Bayes} recursively in time. Note that both the prior and the posterior filtering densities in Eq.~\eqref{Kolmogorov} and Eq.~\eqref{Bayes} are defined as the continuum level, which is not practical. Thus, one important research direction in nonlinear filtering is to study how to accurately approximate the prior and the posterior filtering densities. 

\subsection{The state of the art: particle filters}\label{sec:PF}
Particle filters (PF), which is the state of the art in nonlinear filtering, approximate the filtering densities in Eq.~\eqref{Kolmogorov} and Eq.~\eqref{Bayes} using empirical distributions defined by a set of random samples, referred to as ``particles''. To compare with the proposed score-based filter in Section \ref{sec:SNF}, we briefly recall how particle filters use random samples to approximate the filtering densities. 
At the time instant $t$, we assume that we have a set of $M$ particles, denoted by $\{x_{t,m}\}_{m=1}^M$, that follows the posterior filtering density $P(X_{t} | \mathcal{Y}_t)$. The empirical distribution for approximating $P(X_{t} | \mathcal{Y}_t)$ is given by
\begin{equation}\label{eq:emp}
P(X_t | \mathcal{Y}_t) \approx P^{\rm PF}_{t|t}(X_t) := \f{1}{M}\sum_{m=1}^M \delta_{x_{t,m}}(X_t),
\end{equation}
where $\delta_{x_{t,m}}$ is the Dirac delta function at the $m$-th particle $x_{t,m}$. 
In practice, PF is implemented through the following 3-step procedure to propagate from time $t$ to $t+1$:

\begin{itemize}[leftmargin=15pt]
\item {\em The prediction step.} It is to propagate the particle cloud through the state dynamics and generate a set of predicted particles. For each particle $x_{t,m}$ that represents state $X_t$, we run the state equation in Eq.~\eqref{NLF} to obtain a predicted particle $\tilde{x}_{t+1,m} = f(x_{t,m}, \omega_{t,m})$, where  $\omega_{t,m}$ is a sample of the random variable $\omega_t$.  As a result, we obtain a set of particles  $\{\tilde{x}_{t+1,m}\}_{m=1}^M$ and the corresponding empirical distribution
\begin{equation}\label{PF:Empirical-prediction}
P(X_{t+1}|\mathcal{Y}_t) \approx P^{\rm PF}_{t+1|t}(X_{t+1}) := \f{1}{M}\sum_{m=1}^M \delta_{\tilde{x}_{t+1,m}}(X_{t+1}),
\end{equation}
which approximates the prior filtering density $P(X_{t+1} |\mathcal{Y}_t)$ in Eq.~\eqref{Kolmogorov}.

\item {\em The update step.} It is to incorporates the new observational data $Y_{t+1}$ through Bayesian inference to update the prior filtering density $P_{t+1|t}^{\rm PF}(X_{t+1})$ to the posterior filtering density $\tilde{P}_{t+1|t+1}^{\rm PF}(X_{t+1})$, i.e.,
\begin{equation}\label{PF:Bayes}
P(X_{t+1}|\mathcal{Y}_{t+1}) \approx \tilde{P}^{\rm PF}_{t+1 | t+1}(X_{t+1}) := \sum_{l=1}^L w_{t+1,m} \delta_{\tilde{x}_{t+1,m}}(X_{t+1}),
\end{equation} 
where $\delta_{\tilde{x}_{t+1,m}}$ is the Dirac delta function at the $m$-th predicted particle $\tilde{x}_{t+1,m}$, and the weight $w_{t+1,m} \propto P(M_{t+1}|\tilde{x}_{t+1,m})$ is defined by the likelihood function $P(M_{t+1}|X_{t+1})$ at the $m$-th predicted particle $\tilde{x}_{t+1,m}$.

\item {\em The resampling step.} It is to alleviate the degeneracy issue, in which only a few particles have significant weights while the weights on other particles maybe neglectable. The resampling is often implemented to re-generate a set of equally weighted particles that follows the weighted empirical distribution $\tilde{P}^{\rm PF}_{t+1 | t+1}(X_{t+1})$. We denote the resampled equally weighted particles by $\{x_{t+1,m}\}_{m=1}^M$, which formulate the empirical distribution 
\begin{equation}\label{eq:PFpost}
P(X_{t+1}|\mathcal{Y}_{t+1}) \approx {P}^{\rm PF}_{t+1 | t+1}(X_{t+1}) := \f{1}{M}\sum_{m=1}^M \delta_{{x}_{t+1,m}}(X_{t+1}),
\end{equation}
which is the final approximation of the posterior filtering density at the time instant $t+1$.

\end{itemize}

%
%
%

The challenge of particle filters is the so-called degeneracy issue, especially for high-dimensional nonlinear filtering problems or long-term tracking problems. In these scenarios, the likelihood weights $w_{t+1,m}$ in Eq.~\eqref{PF:Bayes} tend to concentrate on a very small number of particles. As a result, only a few number of particles are taken into account to construct the approximate posterior filtering density in Eq.~\eqref{eq:PFpost}, and the information of the prior filtering density stored in the particles with small weights is gradually ignored. 
The main reason causing the degeneracy issue is the use of a discrete approximation (i.e., the empirical distributions) based on a finite number of pre-chosen samples to characterize the continuous filtering densities in Eq.~\eqref{Kolmogorov} and Eq.~\eqref{Bayes}. Due to the ``curse of dimensionality'', a distribution in a high-dimensional space contains enormous information, and it's very difficult for finite amount of particles to capture the {continuous} characterization of the filtering densities in high-dimensional spaces. This challenge motivated us to exploit recent advances in diffusion models to develop a score-based nonlinear filter that uses a continuous score function to indirectly store the information of the filtering densities. 
%
%

\section{Our method: the score-based filter (SF)}\label{sec:SNF}
This section contains the key components of the proposed method. The score-based diffusion model is briefly introduced in Section \ref{Diffusion-Model}. We introduce the details of the score-based filter in Section \ref{Score-Filter} with the implementation details provided in \ref{app:diff}. 

\subsection{The score-based diffusion model}\label{Diffusion-Model}
The diffusion model is a type of generative machine learning models for generating samples from a target probability density function, denoted by 
\begin{equation}\label{eq:Q}
Q(Z) \;\;\; \text{for}\;\;\; Z \in \mathbb{R}^d,
\end{equation}
where $Z$ is a $d$-dimensional random variable. The key idea is to transform the unknown density $Q(Z)$ to a standard probability distribution, e.g., the standard Gaussian $\mathcal{N}(0, \mathbf{I}_d)$. To this end, a diffusion model first define a forward stochastic differential equation (SDE), i.e., 
\begin{equation}\label{eq:forward}
\hspace{-1cm}\text{\small \bf Forward SDE:}\;\;\;\; Z_\tau = b(\tau) Z_{\tau} d\tau + \sigma(\tau) dW_\tau\;\; \text{ for } \tau \in \mathcal{T} = [0,1],  
\end{equation}
where $\mathcal{T} = [0,1]$ is a pseudo-temporal domain that is different from the real temporal domain in which the nonlinear filtering problem is defined, $W_\tau$ is a standard $d$-dimensional Brownian motion, $b: \mathcal{T} \mapsto \mathbb{R}$ is the drift coefficient,  $\sigma: \mathcal{T} \mapsto \mathbb{R}$ the diffusion coefficient, and the solution $\{Z_\tau\}_{0\leq \tau\leq 1}$ is a diffusion process that takes values in $\mathbb{R}^d$.

The probability density function of the forward process $Z_\tau$ is denoted by
\begin{equation}\label{eq:diffp}
Q_\tau(Z_\tau)\;\; \text{ for }\;\; \tau \in [0,1].
\end{equation}
With a proper definition of $b(\tau)$ and $\sigma(\tau)$ (e.g., see \ref{app:diff}), the forward SDE in Eq.~\eqref{eq:forward} can transform any initial distribution of $Q_0(Z_0)$ to a standard Gaussian variable $Q_1(Z_1) = \mathcal{N}(0, \mathbf{I}_d)$. Therefore, when we set the initial state to be the target random variable, i.e., $Z_0 = Z$ in Eq.~\eqref{eq:forward}, the forward SDE can transform the target distribution $Q(Z)$ to the standard Gaussian distribution $\mathcal{N}(0, \mathbf{I}_d)$. 

It is easy to see that the forward SDE cannot be used to generate samples of the target distribution $Q(Z)$. To do this, the score-based diffusion model shows that there is an equivalent reverse-time SDE to transform the terminal distribution $Q_1(Z_1) = \mathcal{N}(0, \mathbf{I}_d)$ to the initial distribution $Q_0(Z_0)$, i.e.,
\begin{equation}\label{DM:RSDE}
\hspace{-1cm}\text{\small \bf Reverse-time SDE:}\;\;\;\;  d{Z}_\tau = \left[ b(\tau){Z}_\tau - \sigma^2(\tau) S(Z_\tau, \tau)\right] d\tau + \sigma(\tau) d\ola{W}_\tau,
\end{equation}
where $\ola{W}_\tau$ is the backward Brownian motion and $S(Z_\tau, \tau)$ is referred to as the {\em score function}
\begin{equation}\label{eq:exact_score}
\hspace{-0.8cm}\text{\small \bf Score function:}\;\;\;\;  S(Z_{\tau}, \tau) := \nabla_z \log Q_\tau({Z}_\tau),
\end{equation}
that is uniquely determined by the initial distribution $Q_0(Z_0)$ and the coefficients $b(\tau)$, $\sigma(\tau)$.
Note that the SDE in Eq.~\eqref{DM:RSDE} is solved from $\tau=1$ to $\tau=0$. 
In this way, if the score function is given, we can easily generate samples from the target distribution $Q(Z)$ by generating samples from $Q_1(Z_1) = \mathcal{N}(0, \mathbf{I}_d)$ and then solving the reverse-time SDE. 

As such, the problem of generating samples becomes a problem of how to approximate the score function. In practice,
the score function can be estimated by training a score-based model on samples with score matching (Hyvarinen, 2005; Song et al., 2019a). To this end, we train a time-dependent parameterized score-based model, denoted by $\bar{S}(Z_\tau, \tau; {\theta})$, to approximate the exact score function by solving the following optimization problem:
%
%
\begin{equation}\label{optimization:score}
\hat{\theta} = \arg\min_{\theta}\E\Big[   \|S(Z_\tau, \tau) -\bar{S}(Z_\tau, \tau; {\theta}) \|_2^2 \Big],
\end{equation}
where $\theta$ denotes the set of the tuning parameters (e.g., neural network weights) for the approximate score function. The implementation details about the loss function is given in \ref{app:loss}.
Once the approximate score function $\bar{S}(Z_\tau, \tau; {\theta})$ is well trained, it can be substituted into the Reverse-time SDE to generate any number of samples from the target distribution.

\subsection{The methodology of the score-based filter}\label{Score-Filter}
\begin{figure}[h!]
\centering
\includegraphics[width=0.8\textwidth]{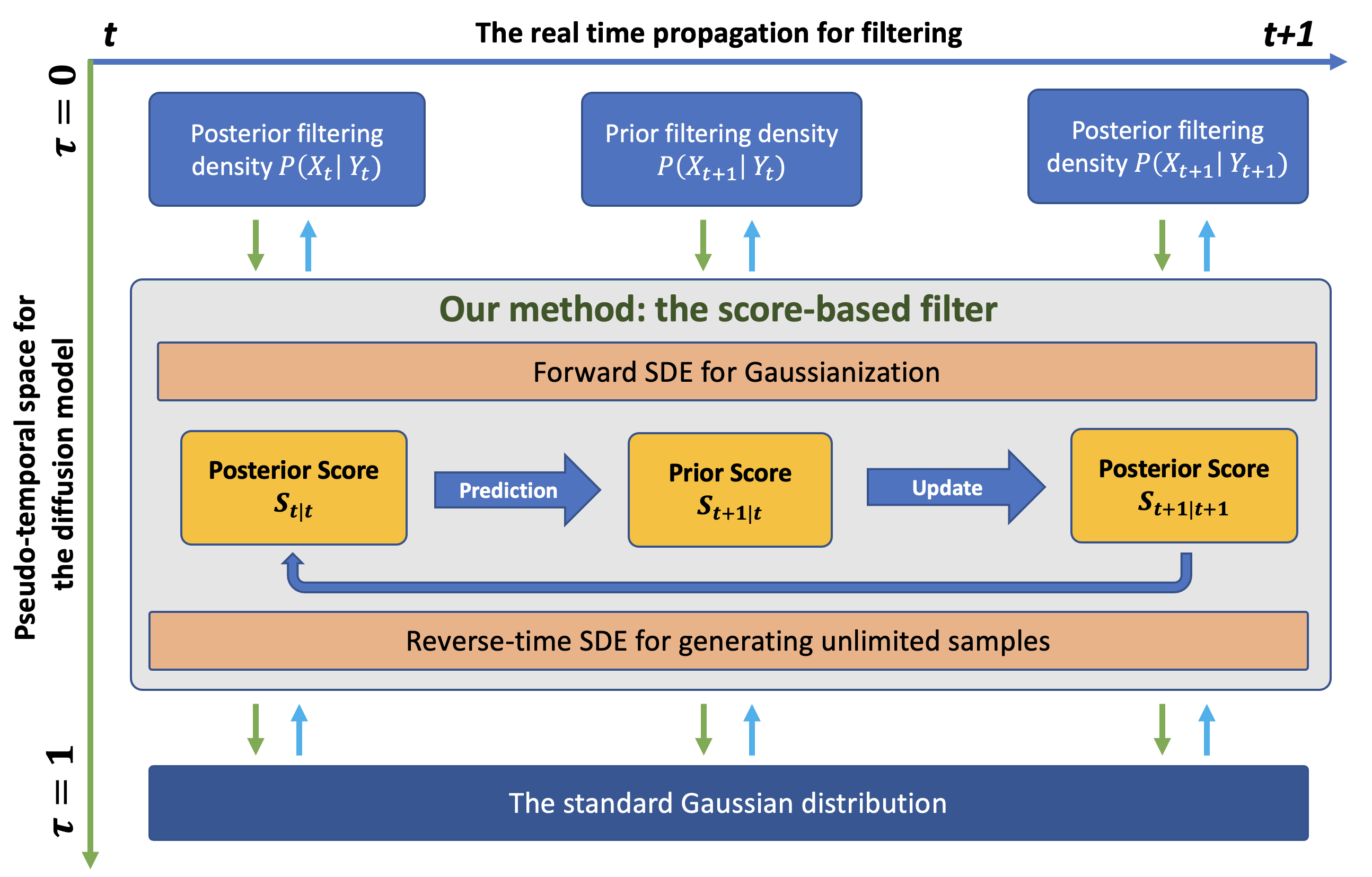}
\caption{The overview of the proposed score-based filter}\label{fig:overview}
\end{figure}

The key idea of the proposed score-based filter is to treat the prior and posterior filtering densities in Eq.~\eqref{Kolmogorov} and Eq.~\eqref{Bayes} as the target distribution $Q(Z)$ in Eq.~\eqref{eq:Q} (or the initial distribution $Q_0(Z_0)$ in Eq.~\eqref{eq:diffp}) in the diffusion model and utilize the score-driven reverse-time SDE to approximate the filtering densities. In other words, we store the information of the filtering densities in the continuous score function.  A systemic overview of the score-based filter is given in Figure \ref{fig:overview}. 

\subsubsection{The relation between the diffusion model and the filtering densities}\label{sec:connect}
Here we discuss how to store the information of the filtering densities in the corresponding score function and how to use the score function to generate unlimited samples of the filtering densities. To proceed, we define the notation 
\begin{equation}\label{eq:Stt}
S_{t|t} (Z_\tau, \tau; \theta)\; \text{ with }\; \tau \in [0,1],
\end{equation}
to represent the exact score function of the posterior filtering density $P(X_t|\mathcal{Y}_t)$ at the time instant $t$. 
%
The diffusion model is related to the filtering density by having the filtering state $X_t$ equal to the initial state $Z_0$ in the forward and reverse-time SDEs in Eq.~\eqref{eq:forward} and Eq.~\eqref{DM:RSDE}, 
\begin{equation}\label{sec:d}
Z_0 = X_t \;\;\xRightarrow{\hspace{0.7cm}}\;\; Q_0(Z_0) = P(X_t | \mathcal{Y}_t),
\end{equation}
where $Q_0(Z_0)$ is the initial distribution of the diffusion model. As shown in \ref{app:diff}, the choice of $b(\tau)$ and $\sigma(\tau)$ in Eq.~\eqref{eq:forward} can ensure that the diffusion model can transform any initial distribution $Q_0$ to the standard Gaussian distribution, 
 we can see that the score function implicitly defines an invertible mapping between the filtering density and the standard Gaussian density, i.e.,
\begin{equation}\label{eq:map}
{\Pi}_{S_{t|t}}\big( \mathcal{N}(0, \mathbf{I}_d) \big) = P(X_{t} | \mathcal{Y}_t)\;\; \text{ and }\;\; 
{\Pi}^{-1}_{S_{t|t}}\big( P(X_{t} | \mathcal{Y}_t) \big) = \mathcal{N}(0, \mathbf{I}_d),
\end{equation}
which indicates that the complete information of the filtering densities can be stored in the score function. This is the key property we will exploit to develop the score-based filter. In Section \ref{sec:sf_pred} and \ref{sec:sf_update}, we will discuss the procedure of dynamically updating the approximate score function $\bar{S}(Z_\tau, \tau; \theta)$.

\subsubsection{The prediction step of the score-based filter}\label{sec:sf_pred}
We intend to evolve the score function $\bar{S}_{t|t}$ associated with $P(X_t|\mathcal{Y}_t)$ to the score function $\bar{S}_{t+1|t}$ associated with the prior filtering density $P(X_{t+1}|\mathcal{Y}_t)$. To achieve this, the prediction step consists of three stages:
\begin{itemize}[leftmargin=15pt]
\item Drawing $J$ samples from $P(X_t | \mathcal{Y}_t)$ by solving the reverse-time SDE in Eq.~\eqref{DM:RSDE} using the score function $\bar{S}_{t|t}$. The samples are denoted by $\{x_{t,j}\}_{j=1}^J$. Unlike the particle filter, we can draw unlimited amount of samples using the diffusion model. 

\item Run the state equation in Eq.~\eqref{NLF} to obtain a predicted samples $\tilde{x}_{t+1,j} = f(x_{t,j}, \omega_{t,j})$, where  $\omega_{t,j}$ is a sample of the random variable $\omega_t$.  

\item Update the score function $\bar{S}_{t|t}$ to $\bar{S}_{t+1|t}$ for the prior filtering density $P(X_{t+1}|\mathcal{Y}_t)$ using the sample set $\{\tilde{x}_{t+1,j}\}_{j=1}^J$ by solving the optimization problem in Eq.~\eqref{optimization:score}.

\end{itemize}

One may notice that scheme $\tilde{x}_{t+1,j} = f(x_{t,j}, \omega_{t,j})$ is similar to the prediction scheme in Eq.~\eqref{PF:Empirical-prediction} in the particle filter, and $\{\tilde{x}_{t+1,j}\}_{j=1}^J$ form a set of samples for the prior filtering density. However, it's important to recall that the score-based filtering stores the information of the target filtering density in the \textit{score function} instead of the \textit{finite} (or discrete) locations of the particles as in the particle filter, and we can generate unlimited number of data samples through the reverse-time SDE as needed to characterize the target distribution. Therefore, the number $J$ in our score-based filter is an arbitrarily chosen number. In this way, as long as the exact score function is well approximated, which is typically obtained through deep learning, we can generate as many samples as needed to pass into the state dynamical model and create a distribution for the predicted filtering density as smooth as we want. On the other hand, once the number of total particles is chosen at the beginning of a particle filter algorithm, the filtering density can be only characterized by the finite locations of those particles.

\subsubsection{The update step of the score-based filter}\label{sec:sf_update}
We intend to update the score function $\bar{S}_{t+1|t}$ corresponding to the prior filtering density to the score function $\bar{S}_{t+1|t+1}$ corresponding to the posterior filtering density. Unlike the prediction step in which we can generate unlimited samples $\{\tilde{x}_{t+1,j}\}_{j=1}^J$ to help us train the score function $\bar{S}_{t+1|t}$, we do not have access to samples from the posterior filtering density. Therefore, we propose to analytically add the likelihood information to the current score $\bar{S}_{t+1|t}$ to define the score $\bar{S}_{t+1|t+1}$ for the posterior filtering density $P(X_{t+1}|\mathcal{Y}_{t+1})$. 

Specifically, we take the gradient of the log likelihood of the posterior filtering density defined in Eq.~\eqref{Bayes} and obtain, 
\begin{equation}\label{eq:log_post}
\nabla_x \log P( X_{t+1}|\mathcal{Y}_{t+1}) =  \nabla_x \log P( X_{t+1}|\mathcal{Y}_{t}) +  \nabla_x \log P(Y_{t+1} | X_{t+1}),
\end{equation}
where the gradient is taken with respect to the state variable at $X_{t+1}$. According to the discussion in Section \ref{sec:connect}, the exact score functions $S_{t+1|t+1}$ and $S_{t+1|t}$ satisfy the following constraints:
\begin{itemize}[leftmargin=50pt]
\item[\bf (C1):]
$
S_{t+1|t}(Z_0, 0) = \nabla_x \log P( X_{t+1}|\mathcal{Y}_{t})\;\; \text{ and }\;\;
S_{t+1|t+1}(Z_0, 0) = \nabla_x \log P( X_{t+1}|\mathcal{Y}_{t+1}),
$
\item[\bf (C2):]
$
 {\Pi}^{-1}_{S_{t+1|t}}\big( P(X_{t+1} | \mathcal{Y}_t) \big) = \mathcal{N}(0, \mathbf{I}_d)\;\; \text{ and }\;\; 
{\Pi}^{-1}_{S_{t+1|t+1}}\big( P(X_{t+1} | \mathcal{Y}_{t+1}) \big) = \mathcal{N}(0, \mathbf{I}_d),
$
\end{itemize}
when setting $Z_0 = X_{t+1}$. We combine the Eq.~\eqref{eq:log_post} and the constraints $\bf (C1)$, $\bf (C2)$ to 
propose an approximation of $S_{t+1|t+1}$ of the form
\begin{equation}\label{eq:d}
\bar{S}_{t+1|t+1}(Z_\tau, \tau; \theta) := \bar{S}_{t+1|t}(Z_\tau, \tau; \theta) + h(\tau) \nabla_z \log P(Y_{t+1} | Z_\tau),
\end{equation}
where $\bar{S}_{t+1|t}(Z_\tau, \tau; \theta)$ is from the prediction step, $\nabla_z \log P(Y_{t+1} | Z_0) = \nabla_x \log P(Y_{t+1} | X_{t+1})$ (if $Z_0 = X_{t+1}$) is analytically defined in Eq.~\eqref{Likelihood}, and $h(\tau)$ is a damping function satisfying
\begin{equation}\label{eq:ht}
h(\tau) \text{ is monotonically decreasing in } [0,1] \text{ with } h(0) = 1 \text{ and } h(1) = 0.
\end{equation}
We use $h(\tau) = 1- \tau$ for $\tau \in [0,1]$ in the numerical examples in Section \ref{sec:exam}. We remark that there are multiple choices of the damping function $h(\tau)$ that satisfying Eq.~\eqref{eq:ht}. How to define the optimal $h(\tau)$ is still an open question that will be considered in our future work.

We observe that the definition of $\bar{S}_{t+1|t+1}(Z_\tau, \tau; \theta)$ is compatible with the constraints $\bf (C1)$, $\bf (C2)$. Intuitively, the information of new observation data in the likelihood function is gradually injected into the diffusion model (or the score function) at the early dynamics (i.e., $\tau$ is small) of the forward SDE during which the deterministic drift (determined by $b(\tau)$) dominates the dynamics. 
When the pseudo-time $\tau$ approaches 1, the diffusion term (determined by $\sigma(\tau)$) becomes dominating, the information in the likelihood is already absorbed into the diffusion model so that $h(1) = 0$ can ensure the final state $Z_1$ still follows the standard Gaussian distribution. 
The updated score function $\bar{S}_{t+1|t+1}(Z_\tau, \tau; \theta)$ can be used as the starting point of the prediction step for the time instant $t+1$. 

\subsubsection{Summary of the score-based filter method}\label{sec:sf_summary}

In Algorithm 1, we use a brief pseudo-algorithm to summarize the score-based filter.

\noindent\makebox[\linewidth]{\rule{\textwidth}{0.5pt}}\\
\vspace{-0.5cm}
\newline {\bf Algorithm 1: the score-based filter}\vspace{-0.2cm} \\ 
\noindent\makebox[\linewidth]{\rule{\textwidth}{0.5pt}}
\vspace{-0.3cm}
\newline1:\, {\bf Input}: the state equation $f(X_t, \omega_t)$, the prior density $P(X_0)$;
\vspace{0.1cm}
\newline2: \,{\bf for} $t = 0, \ldots, $
\vspace{0.1cm}
\newline3: \qquad {\bf if} $t = 0$
\vspace{0.1cm}
\newline4: \qquad \quad Generate $J$ samples $\{x_{0,j}\}_{j=1}^J$ from $P(X_0)$;
\vspace{0.1cm}
\newline5:  \qquad \quad Train the score function $S_{0|0}$ using the sample set $\{x_{0,j}\}_{j=1}^J$;
\vspace{0.1cm}
\newline6: \qquad {\bf else}
\vspace{0.1cm}
\newline7: \qquad \quad\; Generate $J$ samples $\{x_{t,j}\}_{j=1}^J$ of $P(X_t|\mathcal{Y}_t)$ using the score function $\bar{S}_{t|t}$;
\vspace{0.1cm}
\newline8: \qquad Run the state equation in Eq.~\eqref{NLF} to obtain a predicted samples $\{\tilde{x}_{t+1,j}\}_{j=1}^J $;
\vspace{0.1cm}
\newline9: \qquad Train the score function $\bar{S}_{t+1|t}$ using the sample set $\{\tilde{x}_{t+1,j}\}_{j=1}^J$; 
\vspace{0.1cm}
\newline10: \quad\, Update the score function $\bar{S}_{t+1|t}$ to $\bar{S}_{t+1|t+1}$ using Eq.~\eqref{eq:d};
\newline11: {\bf end}\vspace{-0.1cm} \\
\noindent\makebox[\linewidth]{\rule{\textwidth}{0.5pt}}
 \\

As a novel methodology for solving the nonlinear filtering problem, the score-based filter has the following advantages to guarantee its robust accuracy:

\begin{itemize}[leftmargin=25pt]
\item At any recursive stage of the data assimilation procedure, we can substitute the current score function into the reverse-time SDE in Eq.~\eqref{DM:RSDE} to generate \textit{unlimited} samples from the filtering density and compute any statistics of the current state.  

\item The score function modeled by deep neural network (DNN) can take advantage of deep learning, and the DNN learned score is potentially capable to store \textit{complex} information contained in data and state dynamics.

\item The score-based filter is equipped with an \textit{analytical} update step to gradually inject the data information into the score model. This allows the data information to be sufficiently incorporated into the filtering density.  
\end{itemize}

\section{Numerical experiments}\label{sec:exam}

In this section, we demonstrate the performance of our score-based filter by solving three benchmark nonlinear filtering problems. In the first example, we consider a double-well potential problem, which is a 1-dimensional problem with highly nonlinear state dynamics. In the second example, we solve a bearing-tracking problem, in which the state dynamics is linear, but the measurements are nonlinear observational data, i.e. bearings of the state. This is a benchmark example to examine whether a filtering method is suitable for nonlinear problems. The third example that we shall solve in this section is the Lorenz tracking problem, and it is a well-known challenging problem for all the existing filtering methods when the dimension of the problem becomes high.

\subsection{Example 1: Double-well potential}
\begin{figure}[h!]
\begin{center}
\includegraphics[scale = 0.85]{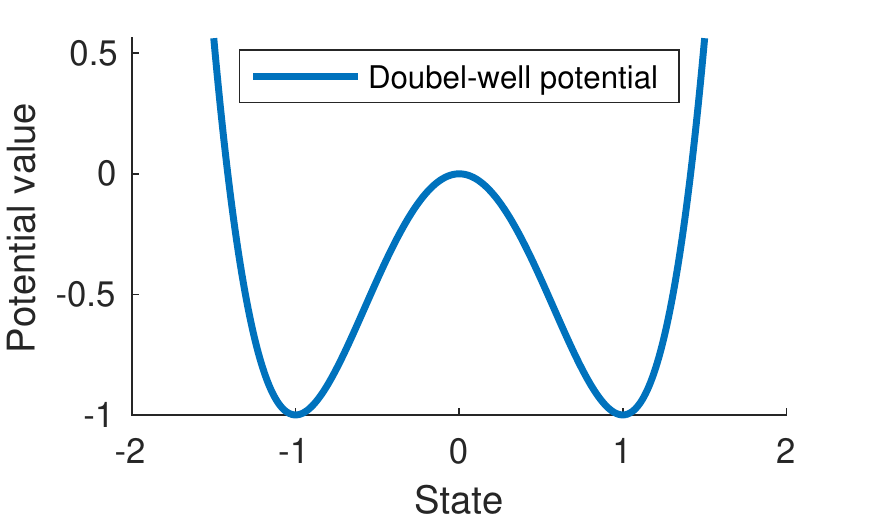}
\end{center}
\caption{Example 1. Double-well potential. }\label{Ex1_Potential_wells} 
\end{figure}

The state dynamics formulated by the double-well potential is given by the following SDE model
\begin{equation}\label{SDE:doubel-well}
dS_t = - 4 S_t (S_t^2 - 1) dt + \beta dB_t,
\end{equation}
where $S$ is the target state, and $B$ is a standard Brownian motion with diffusion coefficient $\beta$. The drift coefficient in Eq. \eqref{SDE:doubel-well} is the derivative of a double-well energy landscape, i.e. $U(x) = (x^4 - 2 x^2)$, which is plotted in Figure \ref{Ex1_Potential_wells}.
The state $S_t$ that follows the dynamics \eqref{SDE:doubel-well} describes a target particle moving on the energy landscape $U$. So there are two stable energy states, i.e. $S_t = 1$ and $S_t = -1$, and there's a force, which is caused by the derivative of the energy potential, that ``drags'' the state towards one of the stable states. 

In this work, we consider the discretized double-well potential model with temporal step-size $\Delta t = 0.1$, and we have the following state process
\begin{equation}\label{Ex1:State}
S_{n+1} =  S_n - 4 \cdot 0.1 \cdot S_n (S_n^2 - 1)  + \beta \sqrt{0.1} \cdot \omega_n.
\end{equation}
In order to track the target particle governed by Eq. \eqref{Ex1:State}, we assume that we have direct observations on $S_{n+1}$, i.e.
\begin{equation}\label{Ex1:Measurement}
M_{n+1} = S_{n+1} + \tau_{n+1},
\end{equation}
where $\tau_{n} \sim N(0, R)$ is the observational noise with standard deviation $R = 0.1$. 

It's easy to track the state while the target stays in the bottom of one of the potential wells. The challenge in solving the nonlinear filtering problem \eqref{Ex1:State} - \eqref{Ex1:Measurement} is that when the target (unexpectedly) switches from one potential well to another, there's a big discrepancy between the state prediction and the measurement data.  In this case, the stability of the optimal filtering algorithm becomes the main issue.

In what follows, we carry out several experiments to compare the performance of our score-based filter with two state-of-the-art optimal filtering methods, i.e. the particle filter and the ensemble Kalmen filter. To implement the score-based filter, we use the sliced score-matching method (see \cite{NEURIPS2019_3001ef25, song2021scorebased}) to solve the diffusion model problem and train the score model with a $50$ neuron - $2$ layer neural network. The sampling procedure through the reverse-time SDE is implemented by the Euler-Maruyama scheme with $K=600$ discretization time steps. Recall that the number of samples that we generate through the reverse-time SDE can be arbitrarily chosen and the computational cost for generating those samples is small. The particle filter that we compare with in this work is the auxiliary particle filter, which is a popular particle filter method with a moderate-cost resampling procedure to improve the performance of the standard bootstrap particle filter, and we use $1000$ particles to produce an empirical approximation for the filtering density. To implement the ensemble Kalman filter, we use $100$ Kalman filter samples in the ensemble. 

In the first numerical experiment, we track the target state over time interval $[0, 10]$, i.e. $100$ time steps, and we let $\beta = 0.3$. 
\begin{figure}[h!]
\centering
\includegraphics[width=0.7\textwidth]{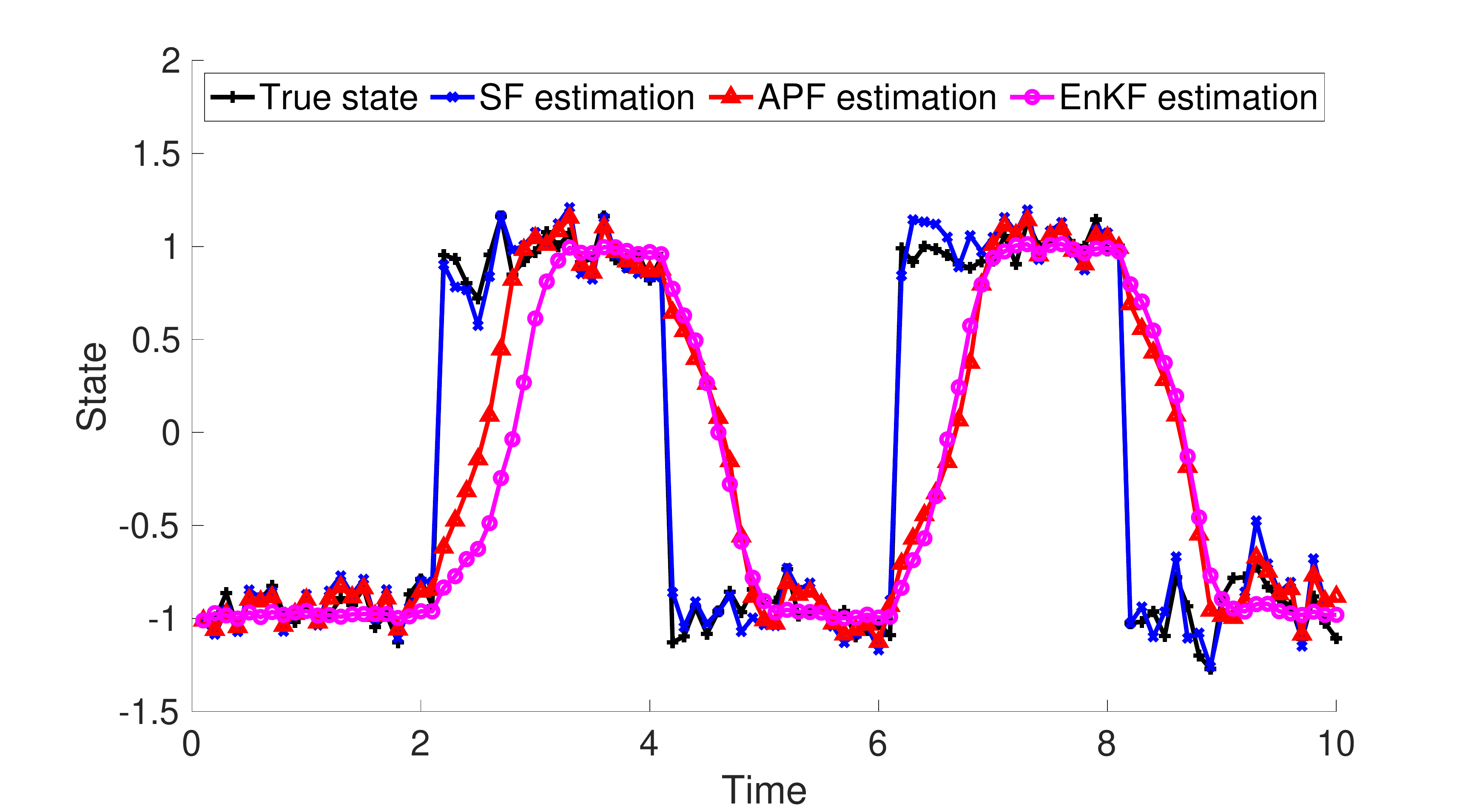}
\vspace{-0.3cm}
\caption{Example 1. Performance comparison with $\beta = 0.3$ }\label{E1_Estimation_Large_Variance} 
\end{figure}
The tracking performance is presented in Figure \ref{E1_Estimation_Large_Variance}, where the black curve (marked by pluses) is the real target state with four unexpected state switches, the blue curve (marked by crosses) is the estimated state obtained by the score-based filter (SF), the red curve (marked by triangles) is the estimated state obtained by the auxiliary particle filter (APF), and the magenta curve (marked by dots) is the estimated state obtained by the ensemble Kalman filter (EnKF). We can see from this figure that all three methods can track the true state while the state stays in the bottom of a potential well. However, when there's a state switch the SF can quickly adjust the changes, and it takes a few estimation steps for the APF and the EnKF to capture the switch. The reason that causes the lagged tracking performance for the APF and the EnKF is that the state model produces strong force that drags the target towards one of the stable energy points, i.e. the bottoms. In this way, the particles (or samples) in the APF (or in the EnKF) will be moved towards the bottom of a potential well. This not only makes the particles (or samples) concentrate at one of the bottoms but also significantly reduces the variance of the predicted filtering density, i.e. the prior. Since both APF and EnKF utilize finite particles (samples) to characterize empirical approximations for the predicted filtering density, narrower distribution bands would result less accurate tail approximations. When the likelihood corresponding to the measurement data lies on the tail of the prior distribution, the data cannot effectively influence the posterior distribution -- due to the poor tail approximation by very limited samples (or no samples at all), which makes the update on the switch much delayed.

\begin{figure}[h!]
\centering
\includegraphics[width=0.7\textwidth]{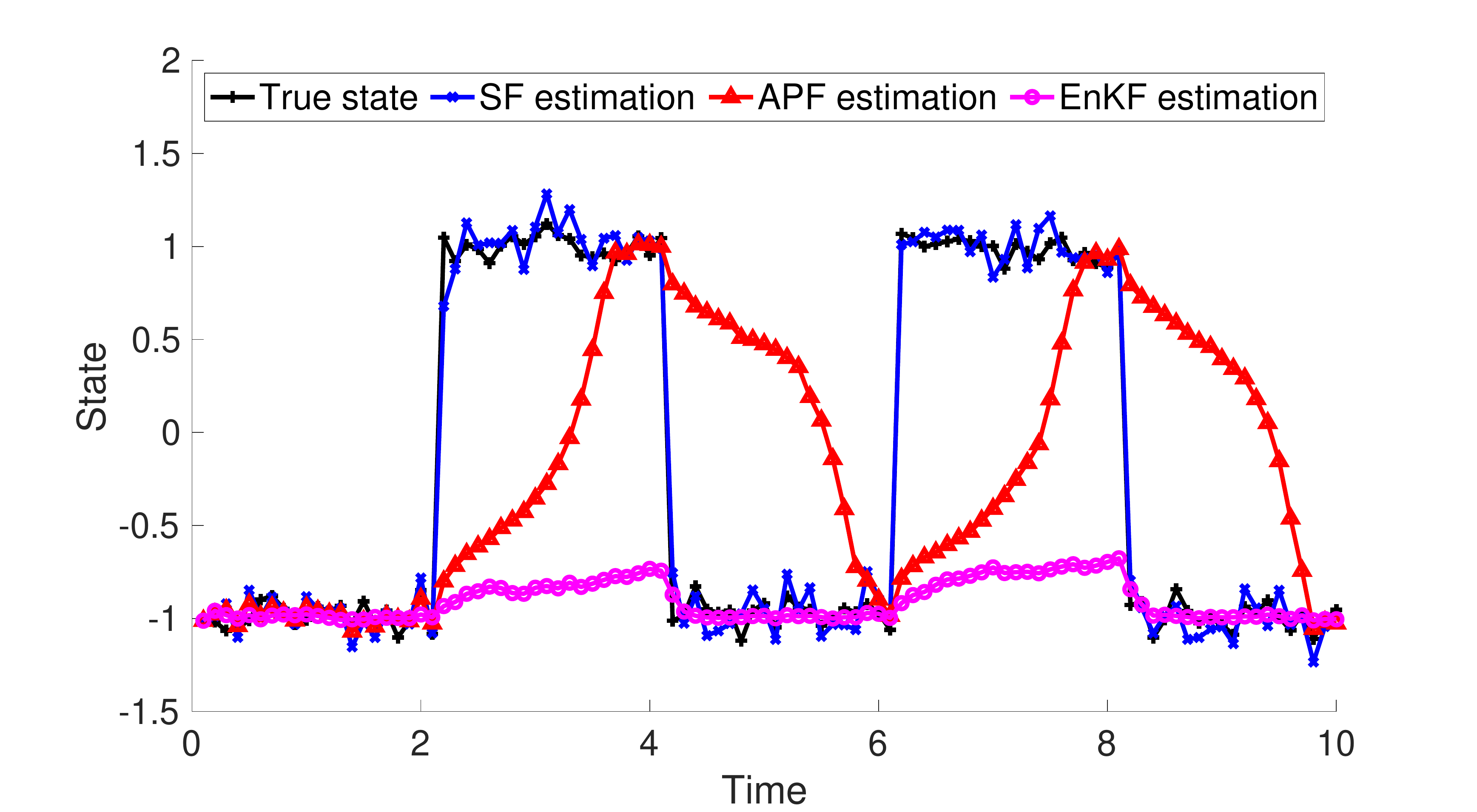}
\vspace{-0.3cm}
\caption{Example 1. Performance comparison with $\beta = 0.2$ }\label{Ex1_Estimation_Small_Variance} 
\end{figure}

To further demonstrate the robust performance of the SF, we reduce the noise variance in the state dynamics to $\beta = 0.2$, and we present the corresponding tracking results of all three methods in Figure \ref{Ex1_Estimation_Small_Variance}.
We can see from this figure that the SF has even more advantageous performance compared with the APF and the EnKF since a smaller noise size, i.e. $\beta = 0.2$, makes the variance of the predicted filtering density obtained by the APF and the EnKF even more narrower. This would result stronger confidence in the predicted state and therefore make the update even less effective. On the other hand, the SF could produce as many samples as needed for the target distribution, which would generate more accurate tail approximations. Moreover, in the SF method the likelihood is continuously added to the posterior score corresponding to the updated filtering density. Therefore, the measurement data information is guaranteed to be incorporated into the updated filtering density through the reverse-time SDE sampler. As a result, the Bayesian update in the SF would be more reliable and more robust.  

\subsection{Example 2: Bearing-only tracking}
In this example, we solve the bearing-only tracking problem, and we consider the following linear dynamical system that models a moving target on the 2-dimensional plane:
\begin{equation}\label{Ex1:State}
S_{n+1} = S_n + A \Delta t + B \sqrt{\Delta t} \cdot \omega_n,
\end{equation}
where $S_n = [x_n, y_n]^{\top}$ describes the position of the target, $A = [v_1, 0; 0, v_2]$ is the velocity matrix that tells how fast the target moves, $B$ is the diffusion coefficient, and $\Delta t$ is the time step. In this example, we let $v_1 = 4$, $v_2 = 6$, $B = [0.2, 0; 0, 0.2]$, and we choose $\Delta t = 0.05$. 

In the bearing-only tracking problem, we can only receive bearings (i.e. angles) for the target state, and we let 
$$M_{n+1} = \arctan \f{y_n - \bar{y}}{x_n - \bar{x}} + \tau_{n+1},$$
where $\tau_{n+1} \sim N(0, R)$ is the observational noise, and $[\bar{x}, \bar{y}]^{\top}$ is an observation platform where we locate the detecter. 
\begin{figure}[h!]
\begin{center}
\includegraphics[scale = 0.7]{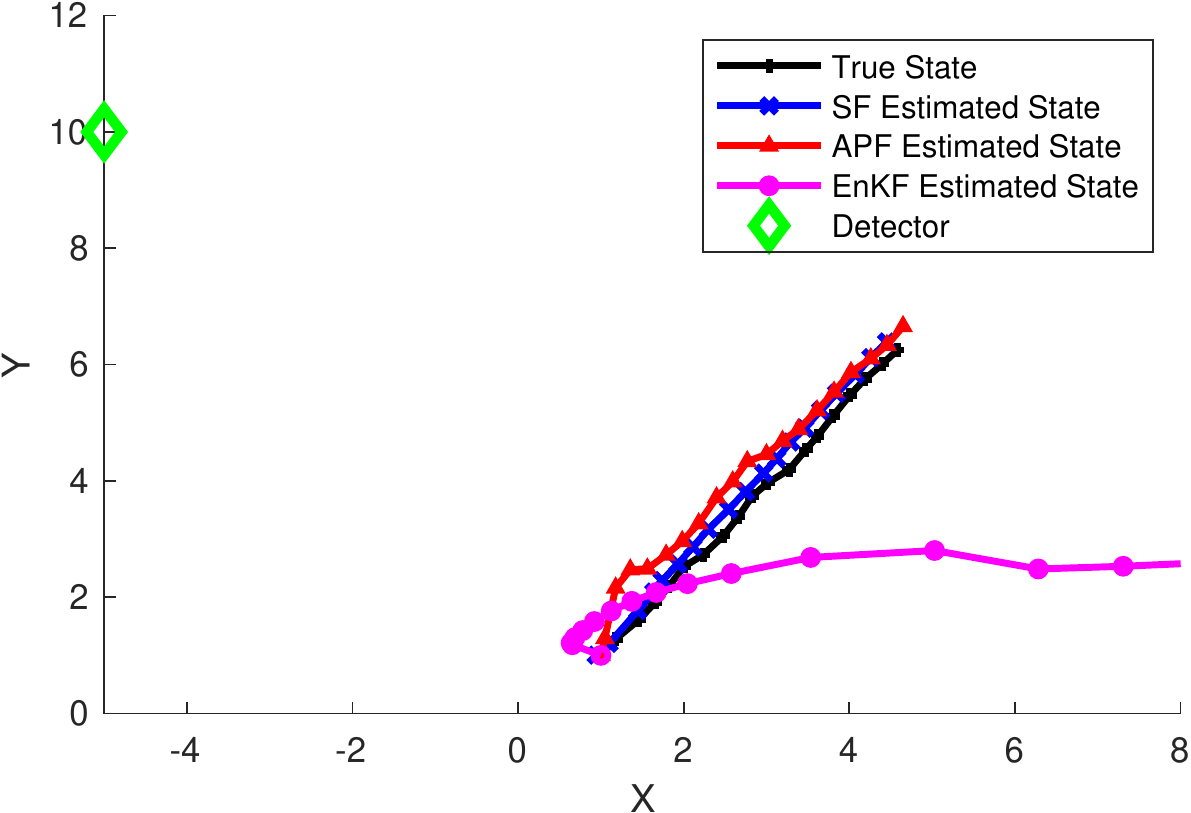}
\end{center}
\caption{Example 2. Comparison of tracking performance }\label{Ex2_Comparison} 
\end{figure}

In Figure \ref{Ex2_Comparison}, we present the tracking performance of the SF, the APF, and the EnKF. All three methods use the same set-up as in Example 1. The green diamond gives the location of the detector platform. The initial position of the target is chosen as $S_0 = [1, 1]^{\top}$. From this figure, we can see that both the SF and the APF can provide accurate estimates for the target state. On the other hand, the EnKF does not work well due to the nonlinear observations \cite{particle-filter}.

\begin{figure}[h!]
\begin{center}
\includegraphics[scale = 1]{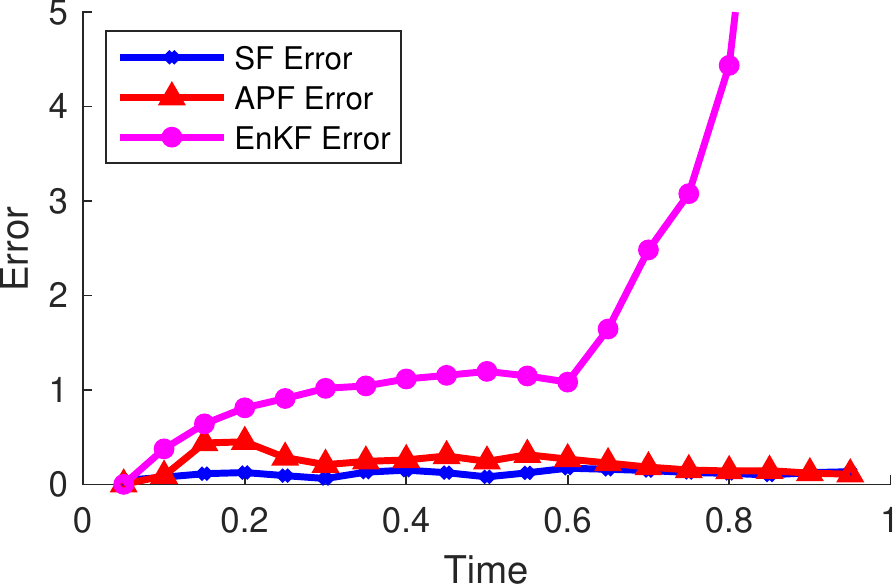}
\end{center}
\caption{Example 2. Comparison of tracking performance }\label{Ex2_Estimation_Errors} 
\end{figure}

To better illustrate the comparison of tracking performance, we plot the tracking errors of the SF, the APF, and the EnKF in Figure \ref{Ex2_Estimation_Errors}.  In this example, we let $R = [0.01, 0; 0, 0.01]$, and the observation platform is chosen at $[\bar{x}, \bar{y}]^{\top} = [-5, 10]^{\top}$, and we track the target for $20$ steps. 

\vspace{0.5em}

From this example, we can see that although the state model is simply a linear dynamical system, the highly nonlinear observational function $\arctan$ makes the Kalman type filters unreliable. This verifies that both the particle filter approach and the score-based filter are suitable for solving the nonlinear filtering problem. 

\subsection{Example 3: Lorenz attractor}

In the third example, we track the state of the Lorenz  96 model described as follows:
\begin{equation}\label{Ex3:Lorenz}
\f{dx_i}{dt} = (x_{i+1} - x_{i-2}) x_{i-1} + F, \qquad i = 1, 2, \cdots, d, \quad d \geq 4,
\end{equation}
where $S_t = [x_1(t), x_2(t), \cdots, x_d(t)]^{\top}$ is a $d$-dimensional target state, and it is assumed that $x_{-1} = x_{d-1}$, $x_{0} = x_{d}$, and $x_{d+1} = x_1$. The term $F$ is a forcing constant. When $F=8$, the Lorenz 96 dynamics \eqref{Ex3:Lorenz} becomes a chaotic system, which makes tracking the state $S_t$ a challenging task for all the existing filtering techniques -- especially when the dimension of the problem is high. 

In this example, we discretize Eq. \eqref{Ex3:Lorenz} through Euler scheme with temporal step-size $\Delta t = 0.01$, and we add a $d$-dimensional white noise (with the standard deviation $0.1$) to perturb the Lorenz 96 model. To make the tracking task more challenging, we assume that our observational data are cubes of the state, i.e.
\begin{equation}\label{Ex3:Measurement}
M_{n+1} = ( S_{n+1} )^{3} + \tau_{n+1}. 
\end{equation}

\begin{figure}[h!] 
\begin{center}
\subfloat[SF \& APF: $x_1$]{\includegraphics[scale = 0.5]{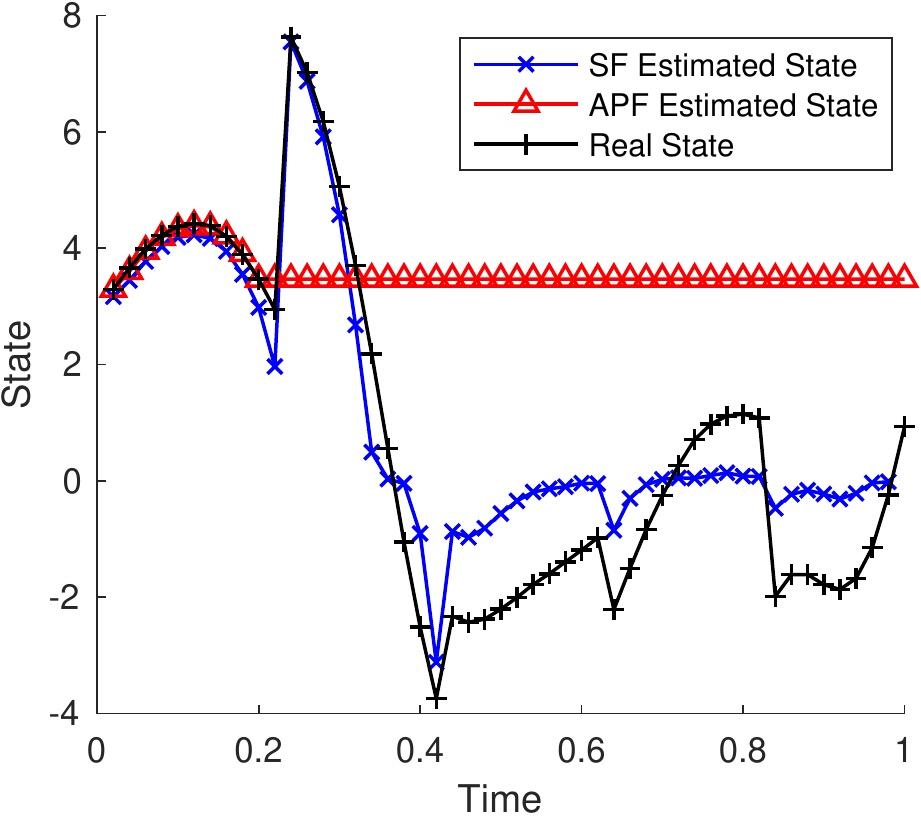} } \quad
\subfloat[SF \& APF: $x_5$]{\includegraphics[scale = 0.5]{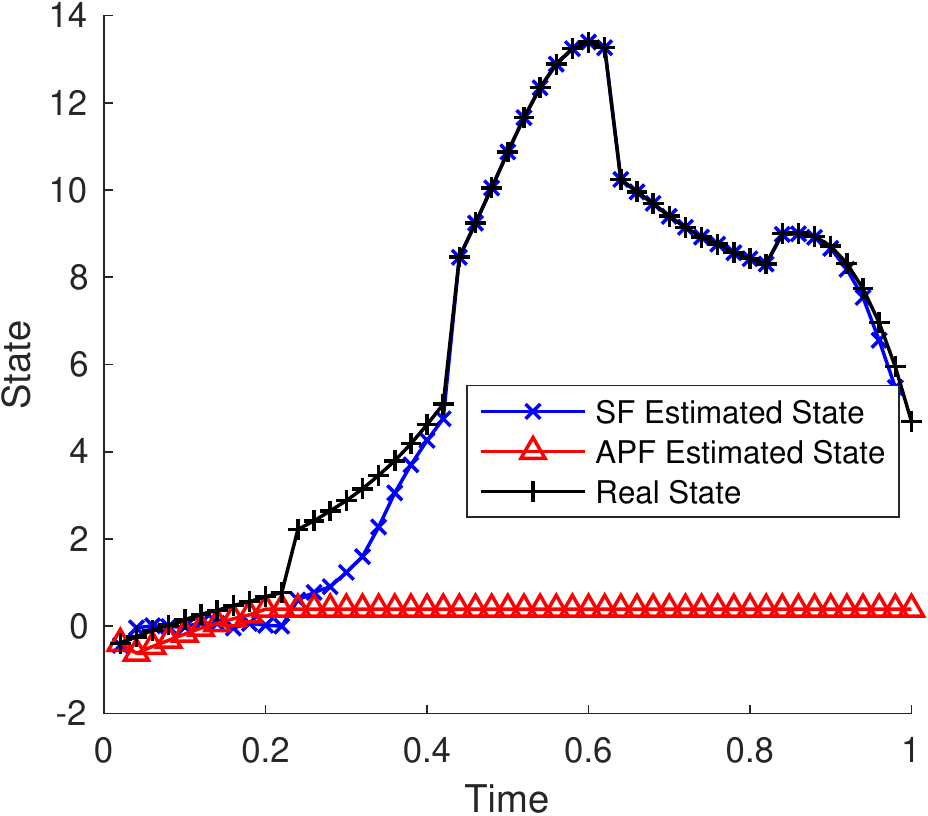} } \quad
\subfloat[SF \& APF: $x_9$]{\includegraphics[scale = 0.5]{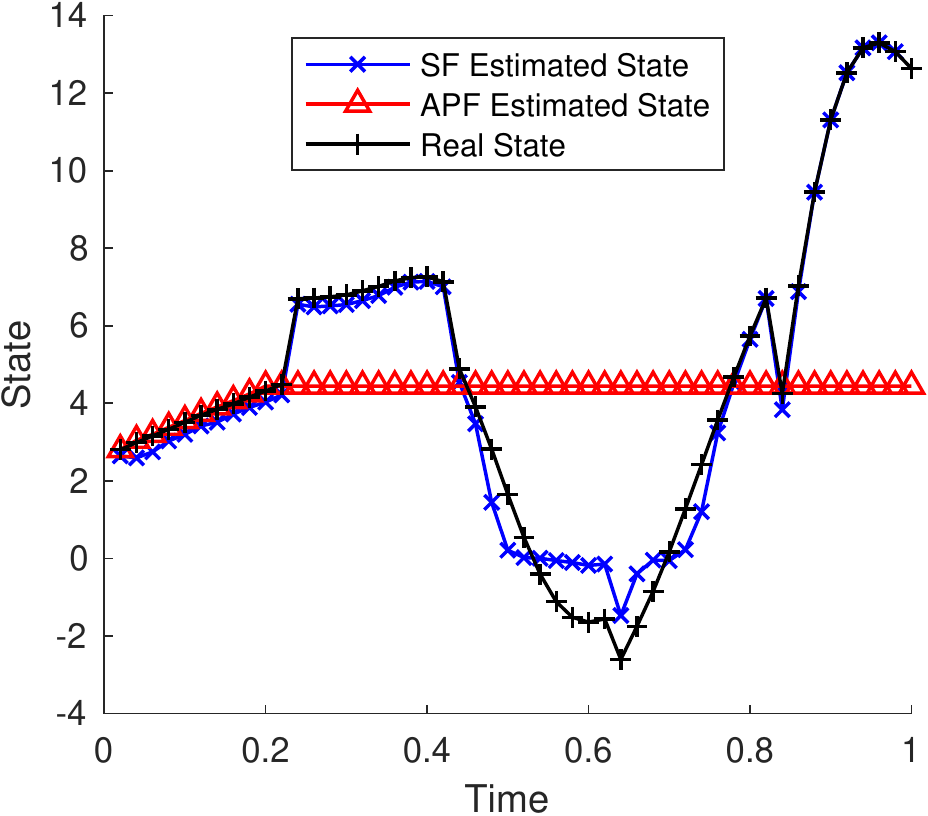} }\\
\subfloat[EnKF: $x_1$]{\includegraphics[scale = 0.5]{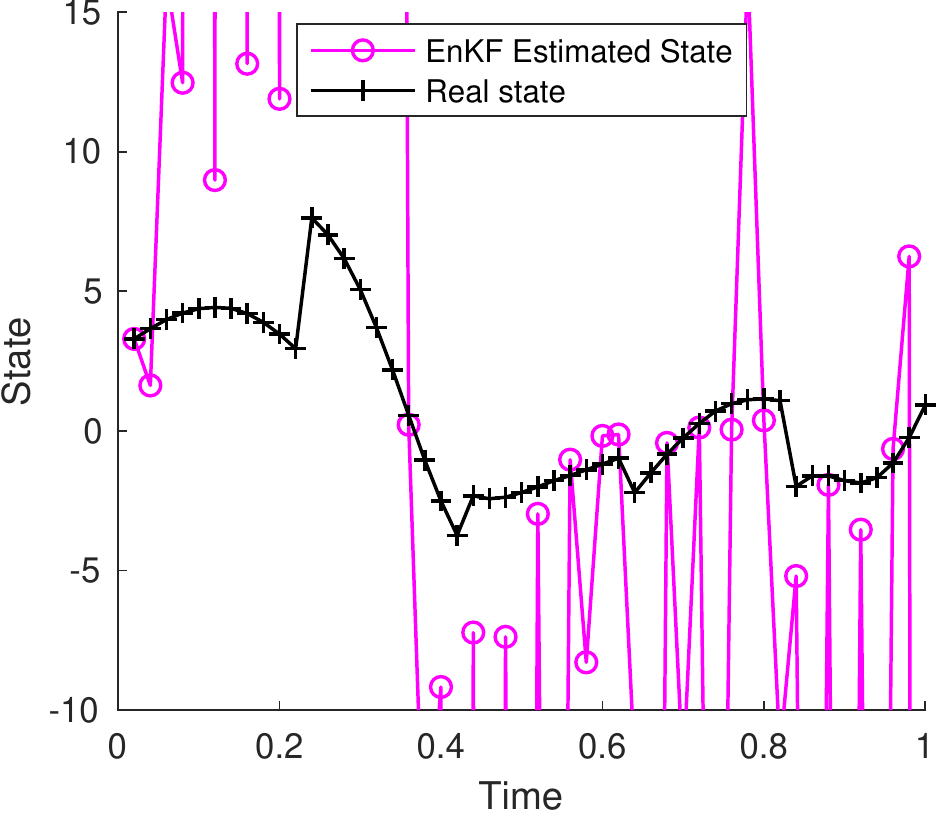} } \quad
\subfloat[EnKF: $x_5$]{\includegraphics[scale = 0.5]{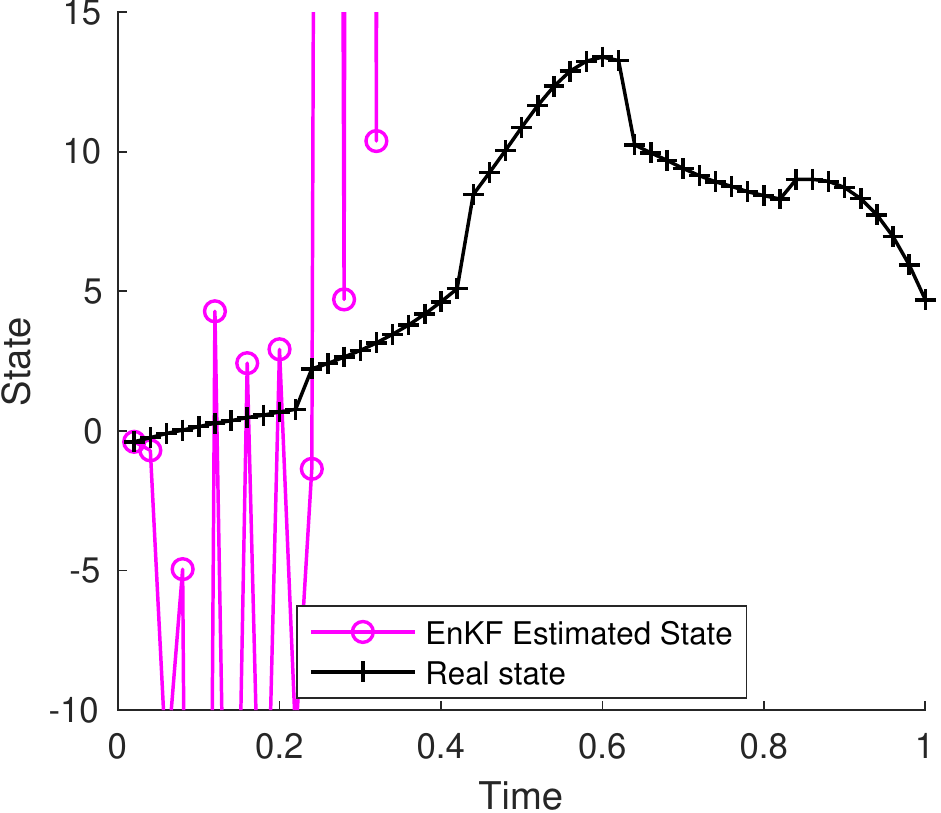} } \quad
\subfloat[EnKF: $x_9$]{\includegraphics[scale = 0.5]{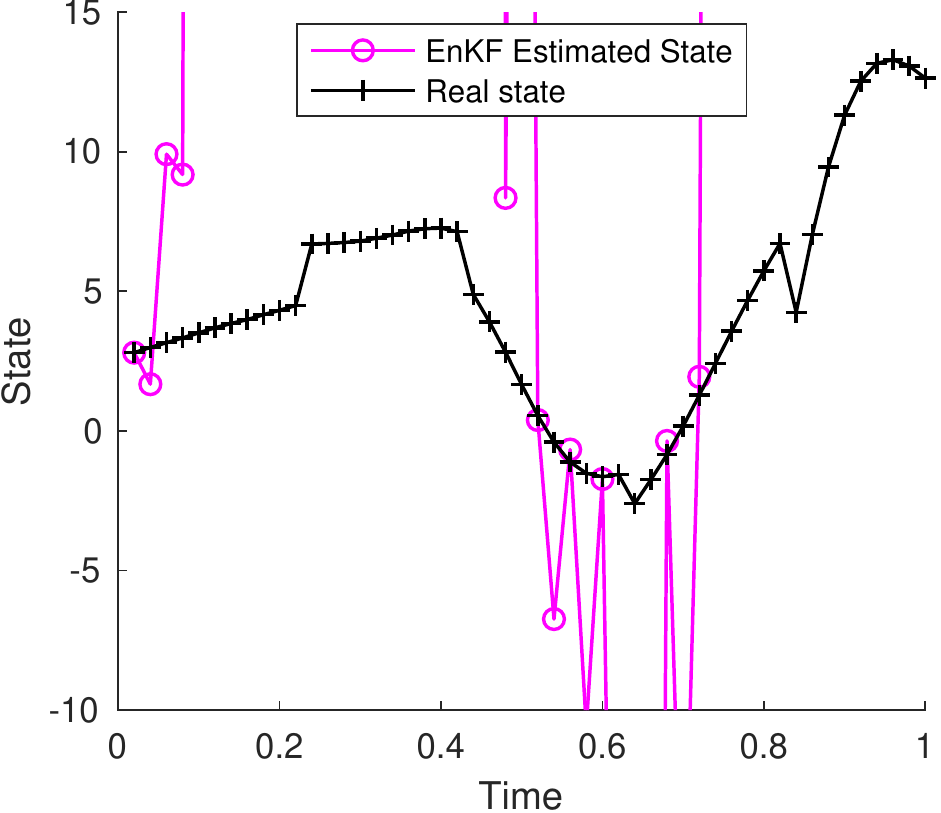} }
\end{center}
\caption{Example 3. State estimation comparison. }\label{10d_Comparison}
\end{figure}

In the first numerical experiment, we track the target state for $100$ time steps in the $10$-dimensional space, i.e. $d = 10$, and we present the state estimation comparison between the SF, the APF, and the EnKF in Figure \ref{10d_Comparison}. Since this is a relatively high dimensional problem, for the SF method we train the score model with a $200$ neuron - $2$ layer neural network, and the sampling procedure through the reverse-time SDE is implemented with $K=800$ discretization time steps. For the APF, we use $20,000$ particles to construct the empirical distribution for the filtering density. For the EnKF,  we use $1000$ Kalman filter samples, which is already a very large number of samples for the EnKF for a 10-dimensional problem.  To make the tracking task even more difficult, we add two random shocks to mimic the severe chaotic behavior of the Lorenz system at time instants $n=21$ and $n=41$, and those unexpected shocks could challenge the stability of different nonlinear filtering methods.

In subplots  (a), (b), and (c) of Figure \ref{10d_Comparison}, we compare the accuracy of state estimation between the SF and the APF in estimating $x_1$, $x_5$, and $x_9$, respectively, where the black curves marked pluses describe the true states, the red curves marked by triangles are the APF estimated states, and the blue curves marked by crosses are the SF estimated states. We can see from those subplots that the APF works well at beginning. However, when the random shock occurs, the large discrepancy between the state prediction and the measurement data makes the particle filter method degenerate. This is indicated by the flat and unresponsive estimation curve. On the other hand, although the SF also suffers from the unexpected random shocks, it

\begin{figure}[h!] 
\begin{center}
\subfloat[State trajectory of real signal]{\includegraphics[scale = 0.6]{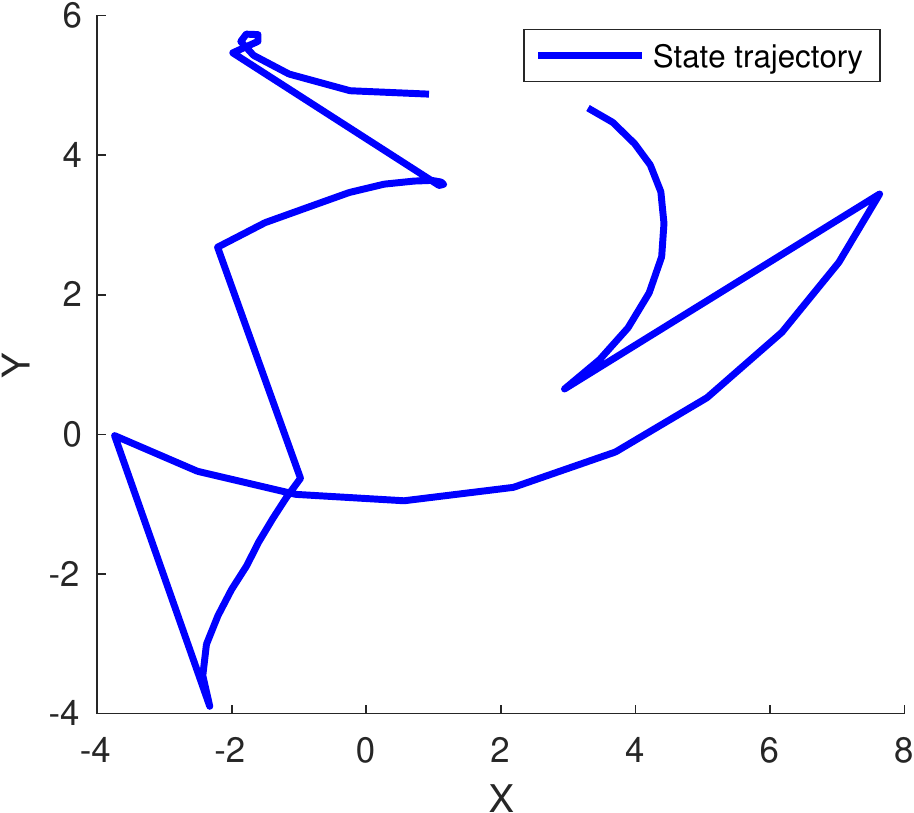} } \qquad
\subfloat[Comparison of state and observation]{\includegraphics[scale = 0.6]{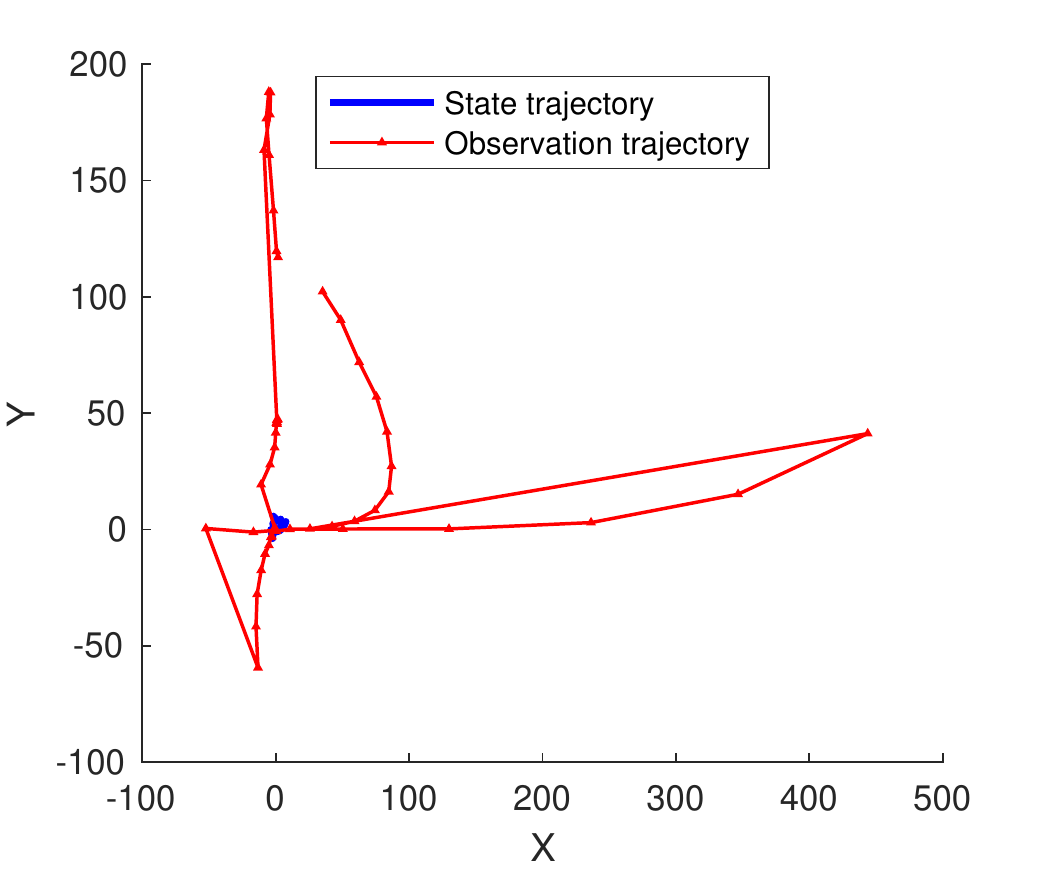} }
\end{center}
\caption{Example 3. Demonstration of state and observation. }\label{State-Observation}
\end{figure}

In Figure \ref{State-Observation}, we plot the first two directions of the original signal $S_t$ in subplot (a), and we compare the cube observations with the real signal in subplot (b), where the observation trajectory is the red curve (marked by triangles) and the true signal is the blue curve. We can see from this figure that the cube measurements are highly nonlinear and only limited data information is contained in the cube observations. Therefore, nonlinear techniques are essential to solve the nonlinear filtering problem \eqref{Ex3:Lorenz}-\eqref{Ex3:Measurement}.

In the second experiment, we further challenge our SF method by solving a $100$-dimensional Lorenz attractor problem.  Note that as a nonlinear filtering method, we need to characterize a non-Gaussian distribution in the $100$-dimensional space, which would certainly encounter the difficulty of ``curse of dimensionality''. To implement the SF, we increase the size of the neural network for the score model to $400$ neurons - $2$ layer, and the sampling procedure through the reverse-time SDE is implemented with $1000$ discretization time steps. For the APF, we utilize $100,000$ particles to approximate the filtering density. For the EnKF, we increase the number of Kalman filter samples to $10,000$. 
\begin{figure}[h!] 
\begin{center}
\subfloat[SF \& APF: $x_{11}$]{\includegraphics[scale = 0.5]{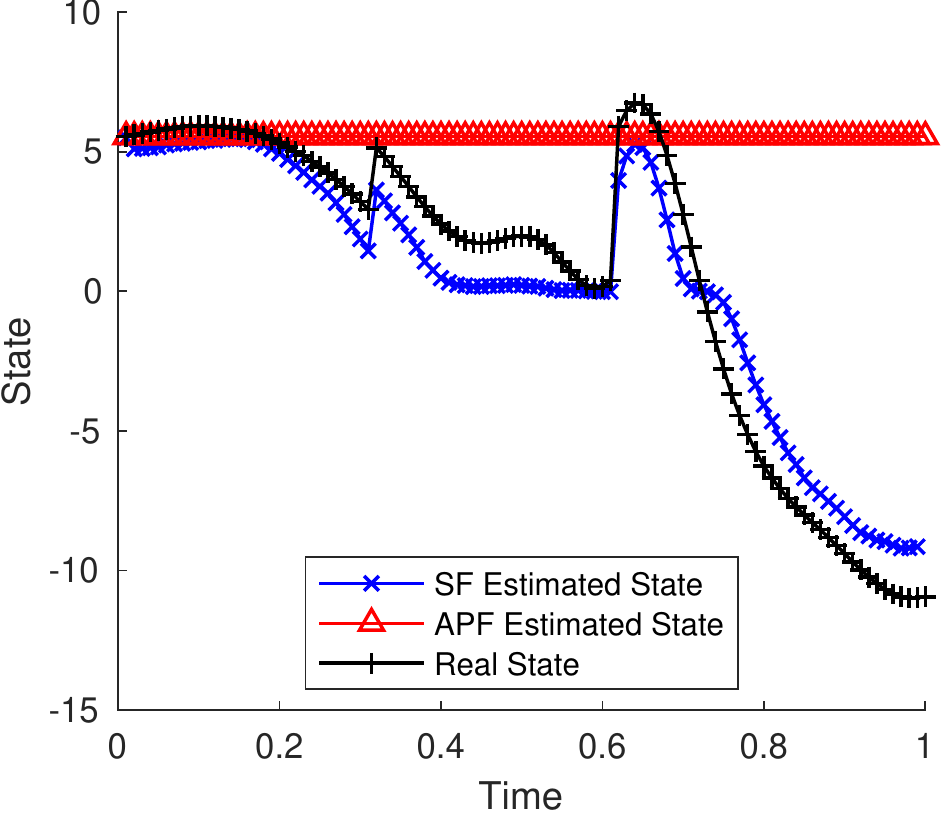} } \quad
\subfloat[SF \& APF: $x_{26}$]{\includegraphics[scale = 0.5]{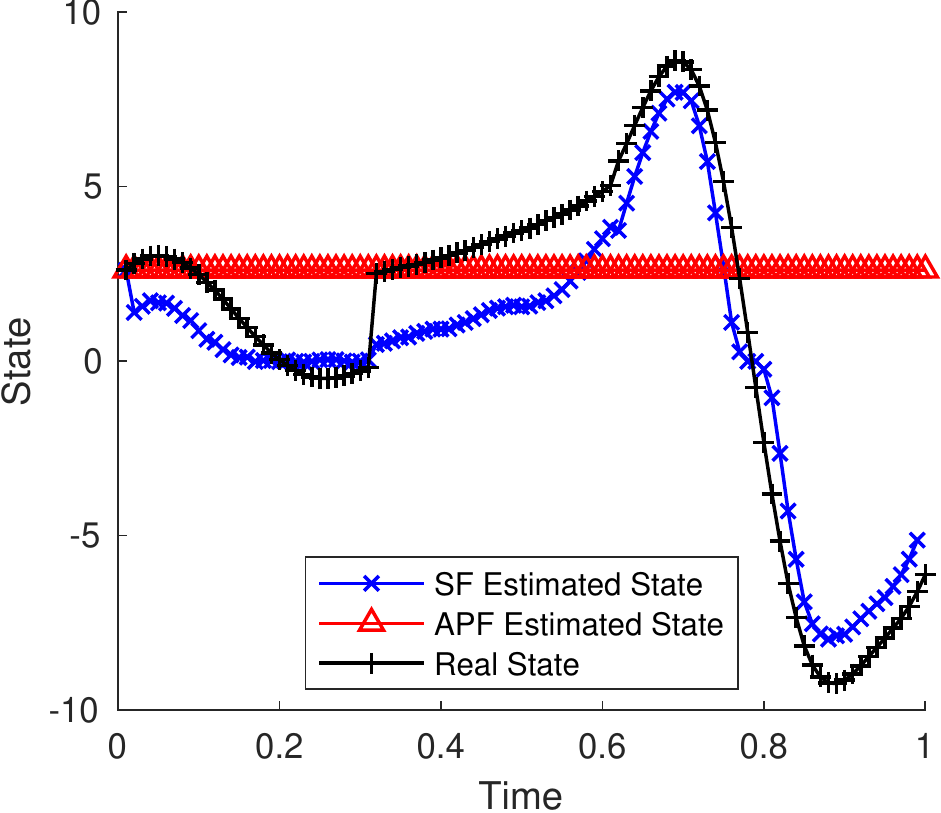} } \quad
\subfloat[SF \& APF: $x_{41}$]{\includegraphics[scale = 0.5]{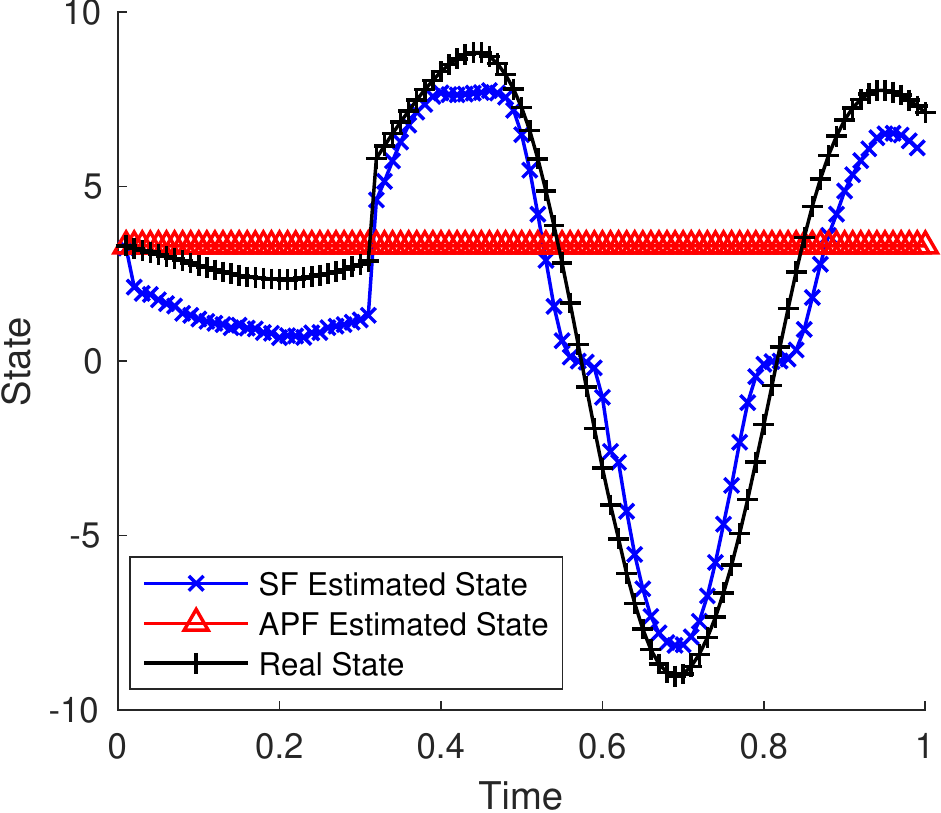} }\\
\subfloat[SF \& APF: $x_{56}$]{\includegraphics[scale = 0.5]{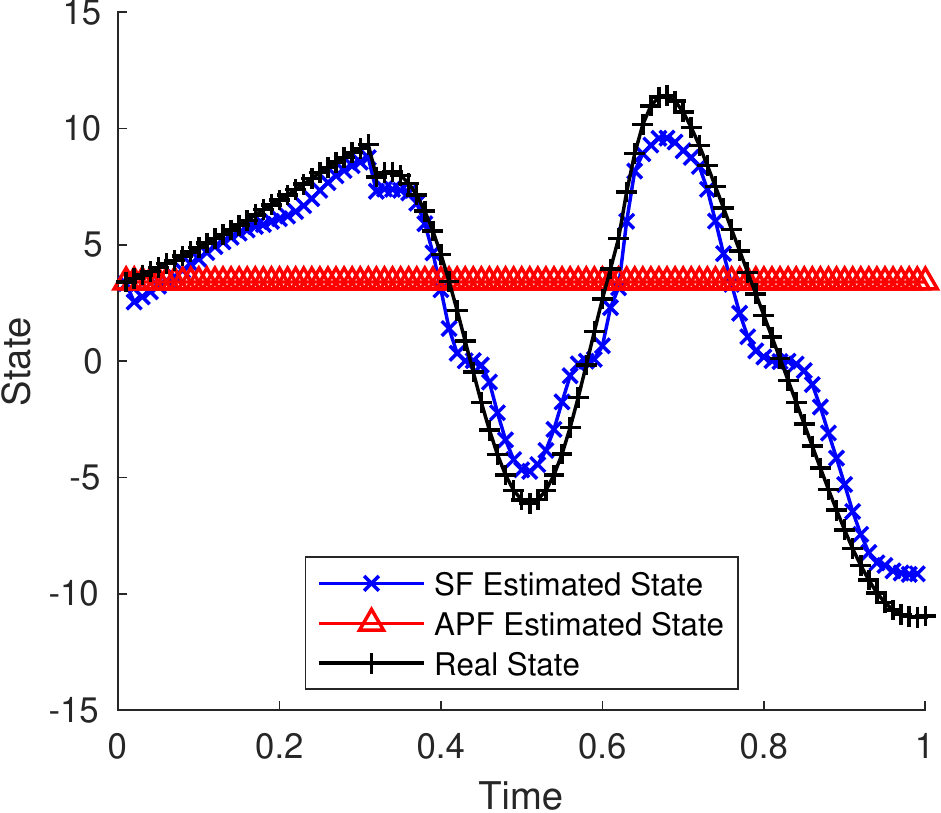} } \quad
\subfloat[SF \& APF: $x_{71}$]{\includegraphics[scale = 0.5]{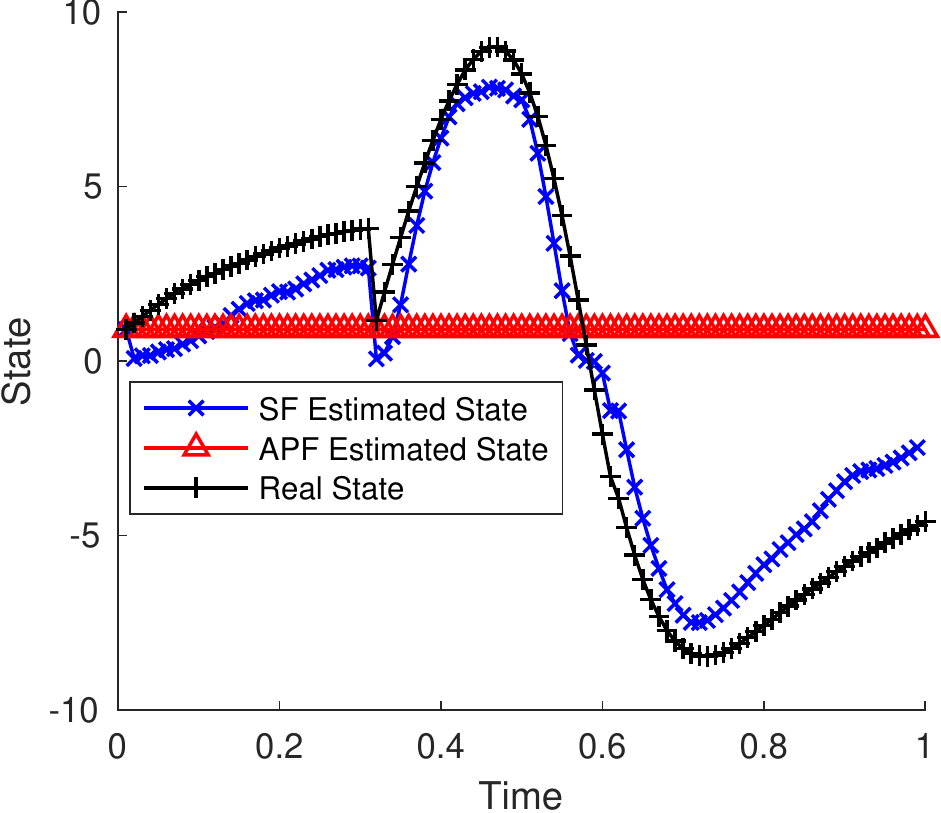} } \quad
\subfloat[SF \& APF: $x_{86}$]{\includegraphics[scale = 0.5]{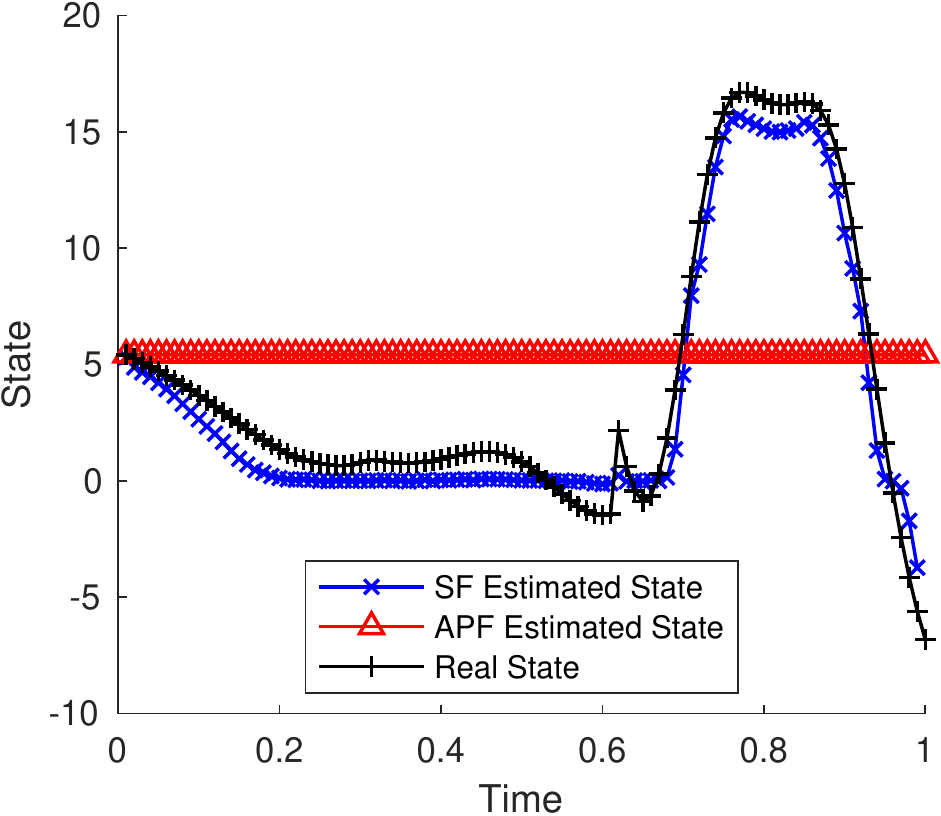} }\\
\end{center}
\caption{Example 3. State estimation comparison: $100$-dimensional case. }\label{100d_Comparison}
\end{figure}
In Figure \ref{100d_Comparison}, we show the comparison results between the SF and the APF in dimensions $x_{11}$, $x_{26}$, $x_{41}$, $x_{56}$, $x_{71}$, and $x_{86}$ in subplots (a) - (f), respectively.   We can see from this figure that the APF suffers from the degeneracy problem from the beginning -- although $100,000$ particles are used to characterize the filtering density, and it fails completely in tracking the $100$-dimensional Lorenz signal. On the other hand, the SF still provides accurate estimates for the target state, and the performance is very stable.
\begin{figure}[h!] 
\begin{center}
\subfloat[EnKF: $x_{11}$]{\includegraphics[scale = 0.5]{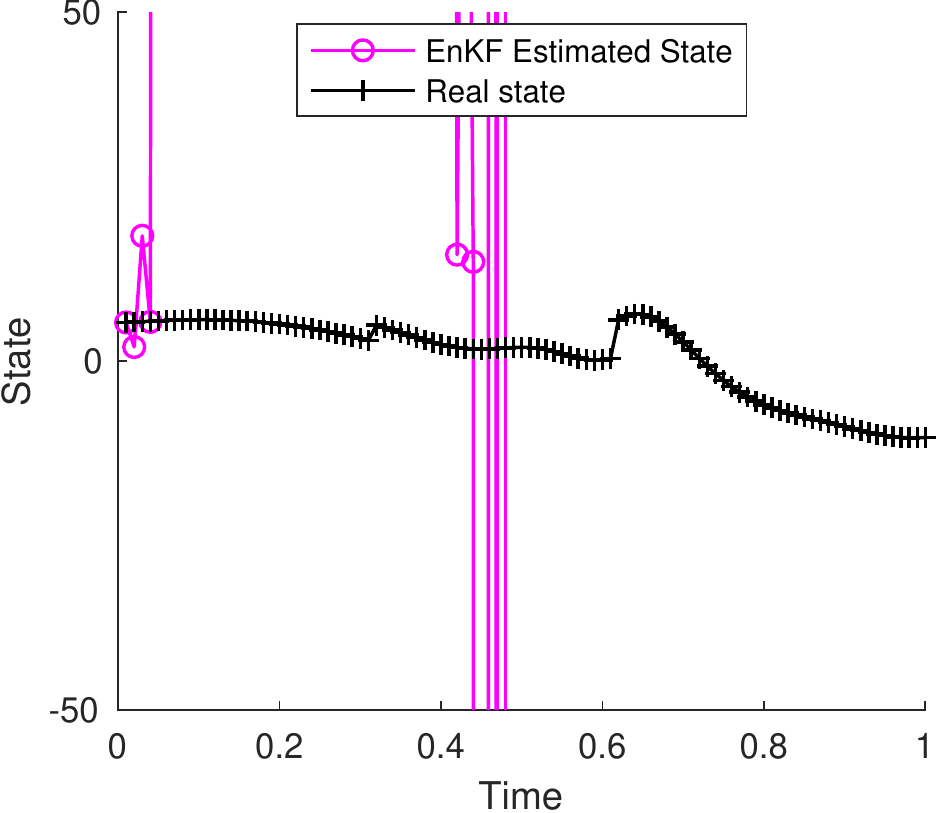} } \quad
\subfloat[EnKF: $x_{26}$]{\includegraphics[scale = 0.5]{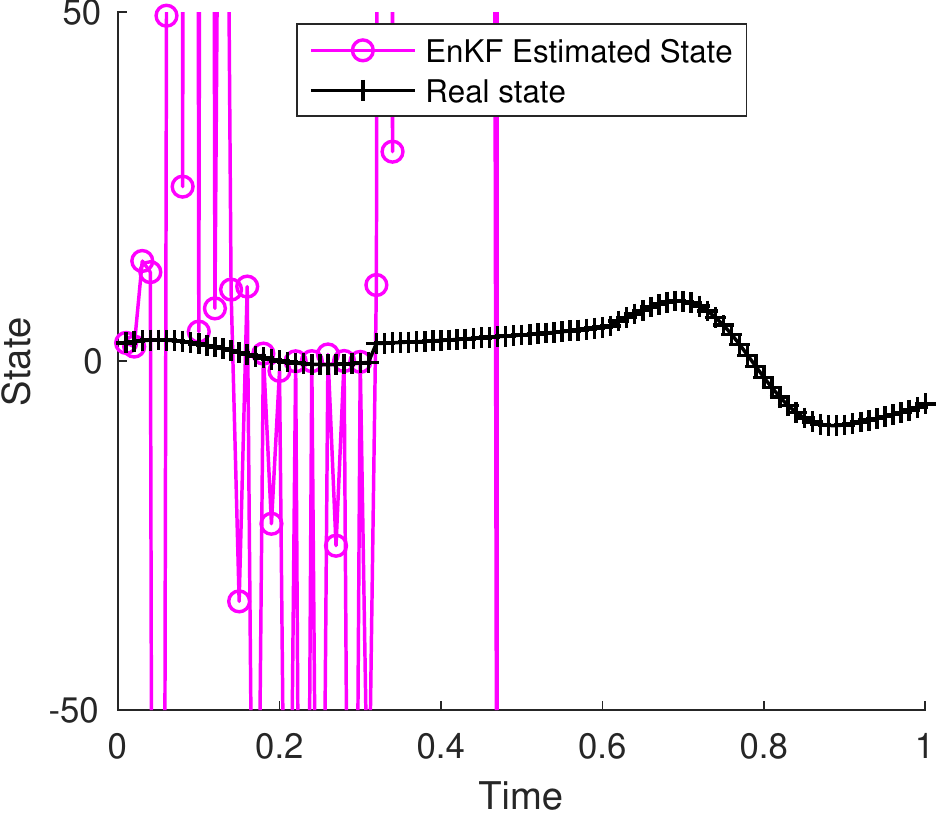} } \quad
\subfloat[EnKF: $x_{41}$]{\includegraphics[scale = 0.5]{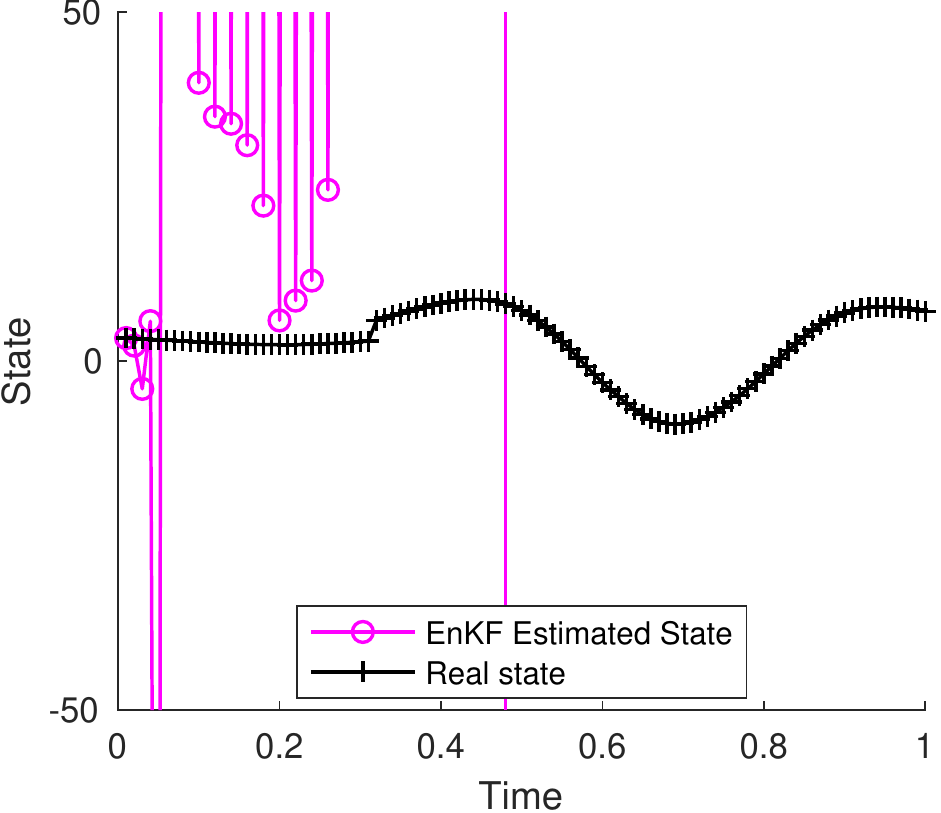} }\\
\subfloat[EnKF: $x_{56}$]{\includegraphics[scale = 0.5]{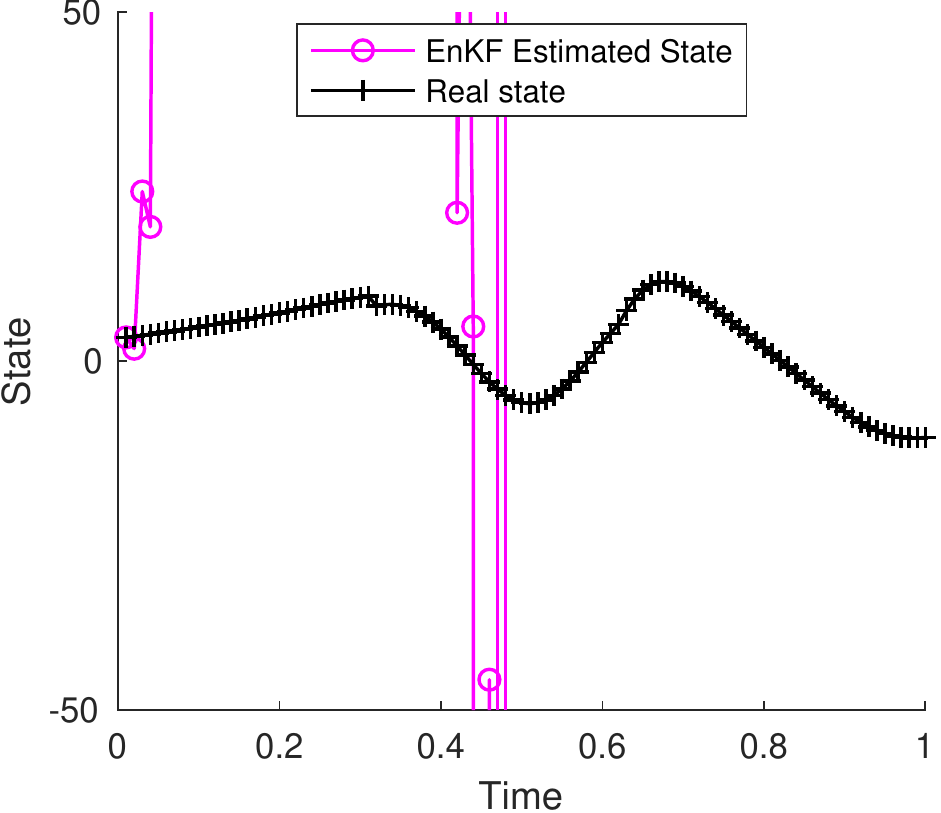} } \quad
\subfloat[EnKF: $x_{71}$]{\includegraphics[scale = 0.5]{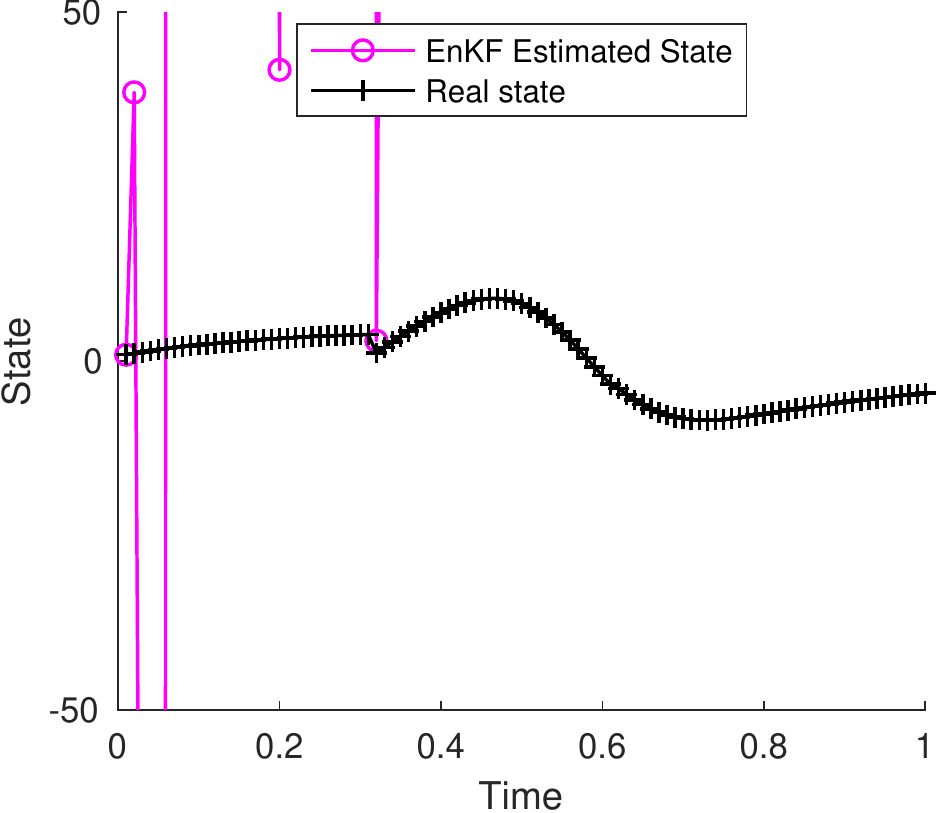} } \quad
\subfloat[EnKF: $x_{86}$]{\includegraphics[scale = 0.5]{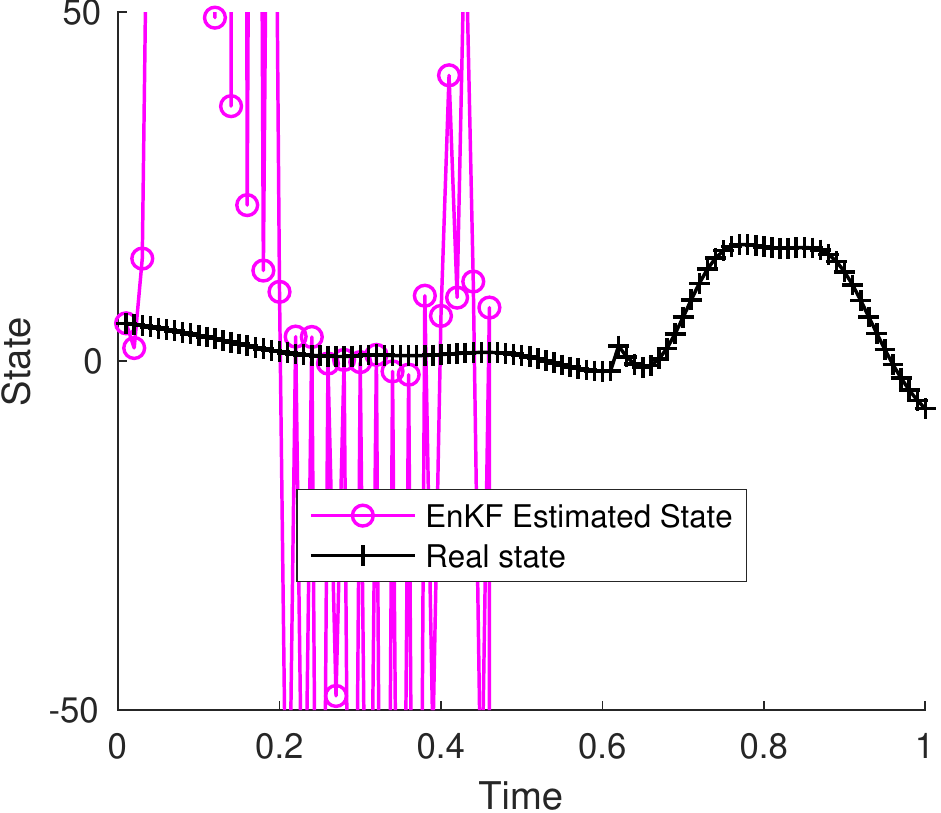} }\\
\end{center}
\caption{Example 3. State estimation by EnKF: $100$-dimensional case. }\label{100d_EnKF}
\end{figure}
In Figure \ref{100d_EnKF}, we present the tracking performance of the EnKF. Similar to the $10$-dimensional tracking experiment, the EnKF does not provide good results for state estimation.

\begin{figure}[h!]
\begin{center}
\includegraphics[scale = 0.8]{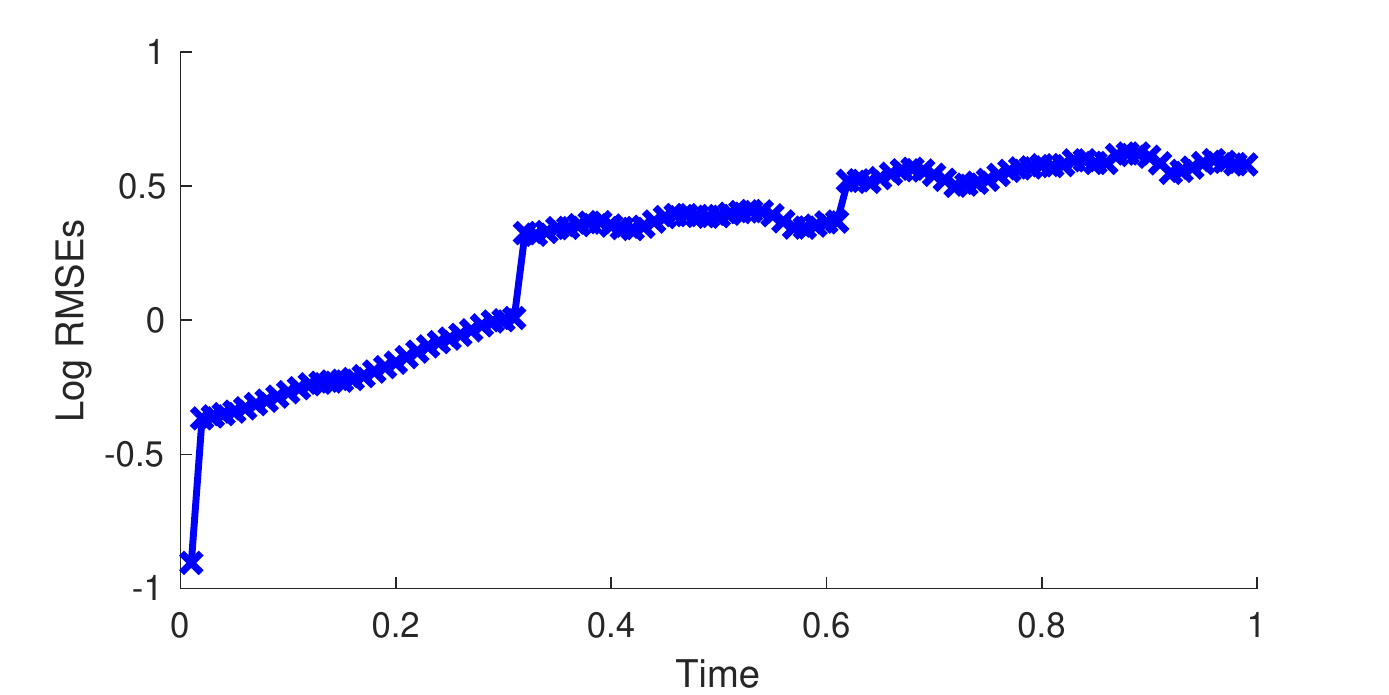}
\end{center}
\caption{Example 3. RMSEs of the SF: $100$-dimensional case. }\label{Ex3_100d_RMSEs} 
\end{figure}

\vspace{0.5em}
To demonstrate the overall performance and the consistent accuracy of the SF in solving the $100$-dimensional Lorenz attractor problem, we repeat the above experiment $20$ times and calculate the root mean square errors (RMSEs). The log of RMSEs with respect to time is presented in Figure \ref{Ex3_100d_RMSEs}, and the error at each time step is the RMSE over all $100$ directions and all $20$ repeated tests. We can see that, although the dimension of the filtering problem is very high and the observations are highly nonlinear, the SF can constantly produce accurate and stable tracking estimates for the chaotic Lorenz 96 system.  

\newpage

\bibliographystyle{plain}
\bibliography{Reference}

\begin{thebibliography}{10}

\bibitem{amit2022segdiff}
Tomer Amit, Tal Shaharbany, Eliya Nachmani, and Lior Wolf.
\newblock Segdiff: Image segmentation with diffusion probabilistic models,
  2022.

\bibitem{MCMC-PF}
C.~Andrieu, A.~Doucet, and R.~Holenstein.
\newblock Particle markov chain monte carlo methods.
\newblock {\em J. R. Statist. Soc. B}, 72(3):269--342, 2010.

\bibitem{DBLP:conf/nips/AustinJHTB21}
Jacob Austin, Daniel~D. Johnson, Jonathan Ho, Daniel Tarlow, and Rianne van~den
  Berg.
\newblock Structured denoising diffusion models in discrete state-spaces.
\newblock In Marc'Aurelio Ranzato, Alina Beygelzimer, Yann~N. Dauphin, Percy
  Liang, and Jennifer~Wortman Vaughan, editors, {\em Advances in Neural
  Information Processing Systems 34: Annual Conference on Neural Information
  Processing Systems 2021, NeurIPS 2021, December 6-14, 2021, virtual}, pages
  17981--17993, 2021.

\bibitem{Bao_CiCP20}
F.~Bao, Y.~Cao, and P.~Maksymovych.
\newblock Backward sde filter for jump diffusion processes and its applications
  in material sciences.
\newblock {\em Communications in Computational Physics}, 27:589--618, 2020.

\bibitem{Bao_AA20}
F.~Bao, Y.~Cao, and J.~Yong.
\newblock Data informed solution estimation for forward backward stochastic
  differential equations.
\newblock {\em Analysis and Applications, to appear}, 2020.

\bibitem{Bao_Cogan20}
F.~Bao, N.~Cogan, A.~Dobreva, and R.~Paus.
\newblock Data assimilation of synthetic data as a novel strategy for
  predicting disease progression in alopecia areata.
\newblock {\em Mathematical Medicine and Biology: A Journal of the IMA}, 2021.

\bibitem{BaoC20142}
Feng Bao, Yanzhao Cao, and Xiaoying Han.
\newblock Forward backward doubly stochastic differential equations and optimal
  filtering of diffusion processes.
\newblock {\em Communications in Mathematical Sciences}, 18(3):635--661, 2020.

\bibitem{Bao_first}
Feng Bao, Yanzhao Cao, Amnon Meir, and Weidong Zhao.
\newblock A first order scheme for backward doubly stochastic differential
  equations.
\newblock {\em SIAM/ASA J. Uncertain. Quantif.}, 4(1):413--445, 2016.

\bibitem{Bao_Zakaid_2015}
Feng Bao, Yanzhao Cao, Clayton Webster, and Guannan Zhang.
\newblock A hybrid sparse-grid approach for nonlinear filtering problems based
  on adaptive-domain of the {Z}akai equation approximations.
\newblock {\em SIAM/ASA J. Uncertain. Quantif.}, 2(1):784--804, 2014.

\bibitem{BSDE_filter}
Feng Bao and Vasileios Maroulas.
\newblock Adaptive meshfree backward {SDE} filter.
\newblock {\em SIAM J. Sci. Comput.}, 39(6):A2664--A2683, 2017.

\bibitem{baranchuk2022labelefficient}
Dmitry Baranchuk, Andrey Voynov, Ivan Rubachev, Valentin Khrulkov, and Artem
  Babenko.
\newblock Label-efficient semantic segmentation with diffusion models.
\newblock In {\em International Conference on Learning Representations}, 2022.

\bibitem{DBLP:conf/cvpr/BrempongKCPM022}
Emmanuel~Asiedu Brempong, Simon Kornblith, Ting Chen, Niki Parmar, Matthias
  Minderer, and Mohammad Norouzi.
\newblock Denoising pretraining for semantic segmentation.
\newblock In {\em {IEEE/CVF} Conference on Computer Vision and Pattern
  Recognition Workshops, {CVPR} Workshops 2022, New Orleans, LA, USA, June
  19-20, 2022}, pages 4174--4185. {IEEE}, 2022.

\bibitem{Multi-target_2007}
Monica~F. Bugallo, Ting Lu, and Petar~M. Djuric.
\newblock Target tracking by multiple particle filtering.
\newblock In {\em 2007 IEEE Aerospace Conference}, pages 1--7, 2007.

\bibitem{DBLP:conf/eccv/CaiYAHBSH20}
Ruojin Cai, Guandao Yang, Hadar Averbuch{-}Elor, Zekun Hao, Serge~J. Belongie,
  Noah Snavely, and Bharath Hariharan.
\newblock Learning gradient fields for shape generation.
\newblock In Andrea Vedaldi, Horst Bischof, Thomas Brox, and Jan{-}Michael
  Frahm, editors, {\em Computer Vision - {ECCV} 2020 - 16th European
  Conference, Glasgow, UK, August 23-28, 2020, Proceedings, Part {III}}, volume
  12348 of {\em Lecture Notes in Computer Science}, pages 364--381. Springer,
  2020.

\bibitem{CT1}
A.~J. Chorin and X.~Tu.
\newblock Implicit sampling for particle filters.
\newblock {\em Proc. Nat. Acad. Sc. USA}, 106:17249--17254, 2009.

\bibitem{NEURIPS2021_49ad23d1}
Prafulla Dhariwal and Alexander Nichol.
\newblock Diffusion models beat gans on image synthesis.
\newblock In M.~Ranzato, A.~Beygelzimer, Y.~Dauphin, P.S. Liang, and J.~Wortman
  Vaughan, editors, {\em Advances in Neural Information Processing Systems},
  volume~34, pages 8780--8794. Curran Associates, Inc., 2021.

\bibitem{Bao_Atom20}
O.~Dyck, M.~Ziatdinov, S.~Jesse, F.~Bao, A.~Yousefzadi~Nobakht, A.~Maksov, B.G.
  Sumpter, R.~Archibald, K.J.H. Law, and S.V. Kalinin.
\newblock Probing potential energy landscapes via electron-beam-induced single
  atom dynamics.
\newblock {\em Acta Materialia}, 203:116508, 2021.

\bibitem{Evense_EnKF}
G.~Evensen.
\newblock The ensemble {K}alman filter for combined state and parameter
  estimation: {M}onte {C}arlo techniques for data assimilation in large
  systems.
\newblock {\em IEEE Control Syst. Mag.}, 29(3):83--104, 2009.

\bibitem{particle-filter}
N.J Gordon, D.J Salmond, and A.F.M. Smith.
\newblock Novel approach to nonlinear/non-gaussian bayesian state estimation.
\newblock {\em IEE PROCEEDING-F}, 140(2):107--113, 1993.

\bibitem{DBLP:journals/corr/abs-2206-09012}
Alexandros Graikos, Nikolay Malkin, Nebojsa Jojic, and Dimitris Samaras.
\newblock Diffusion models as plug-and-play priors.
\newblock {\em CoRR}, abs/2206.09012, 2022.

\bibitem{NEURIPS2020_4c5bcfec}
Jonathan Ho, Ajay Jain, and Pieter Abbeel.
\newblock Denoising diffusion probabilistic models.
\newblock In H.~Larochelle, M.~Ranzato, R.~Hadsell, M.F. Balcan, and H.~Lin,
  editors, {\em Advances in Neural Information Processing Systems}, volume~33,
  pages 6840--6851. Curran Associates, Inc., 2020.

\bibitem{DBLP:journals/jmlr/HoSCFNS22}
Jonathan Ho, Chitwan Saharia, William Chan, David~J. Fleet, Mohammad Norouzi,
  and Tim Salimans.
\newblock Cascaded diffusion models for high fidelity image generation.
\newblock {\em J. Mach. Learn. Res.}, 23:47:1--47:33, 2022.

\bibitem{DBLP:conf/nips/HoogeboomNJFW21}
Emiel Hoogeboom, Didrik Nielsen, Priyank Jaini, Patrick Forr{\'{e}}, and Max
  Welling.
\newblock Argmax flows and multinomial diffusion: Learning categorical
  distributions.
\newblock In Marc'Aurelio Ranzato, Alina Beygelzimer, Yann~N. Dauphin, Percy
  Liang, and Jennifer~Wortman Vaughan, editors, {\em Advances in Neural
  Information Processing Systems 34: Annual Conference on Neural Information
  Processing Systems 2021, NeurIPS 2021, December 6-14, 2021, virtual}, pages
  12454--12465, 2021.

\bibitem{Kang-PF}
Kai Kang, Vasileios Maroulas, Ioannis Schizas, and Feng Bao.
\newblock Improved distributed particle filters for tracking in a wireless
  sensor network.
\newblock {\em Comput. Statist. Data Anal.}, 117:90--108, 2018.

\bibitem{DBLP:conf/iccvw/KawarVE21}
Bahjat Kawar, Gregory Vaksman, and Michael Elad.
\newblock Stochastic image denoising by sampling from the posterior
  distribution.
\newblock In {\em {IEEE/CVF} International Conference on Computer Vision
  Workshops, {ICCVW} 2021, Montreal, BC, Canada, October 11-17, 2021}, pages
  1866--1875. {IEEE}, 2021.

\bibitem{DBLP:journals/corr/abs-2112-05149}
Boah Kim, Inhwa Han, and Jong~Chul Ye.
\newblock Diffusemorph: Unsupervised deformable image registration along
  continuous trajectory using diffusion models.
\newblock {\em CoRR}, abs/2112.05149, 2021.

\bibitem{DBLP:journals/ijon/LiYCCFXLC22}
Haoying Li, Yifan Yang, Meng Chang, Shiqi Chen, Huajun Feng, Zhihai Xu, Qi~Li,
  and Yueting Chen.
\newblock Srdiff: Single image super-resolution with diffusion probabilistic
  models.
\newblock {\em Neurocomputing}, 479:47--59, 2022.

\bibitem{DBLP:journals/corr/abs-2205-14217}
Xiang~Lisa Li, John Thickstun, Ishaan Gulrajani, Percy Liang, and Tatsunori~B.
  Hashimoto.
\newblock Diffusion-lm improves controllable text generation.
\newblock {\em CoRR}, abs/2205.14217, 2022.

\bibitem{DBLP:conf/iccv/LuoH21}
Shitong Luo and Wei Hu.
\newblock Score-based point cloud denoising.
\newblock In {\em 2021 {IEEE/CVF} International Conference on Computer Vision,
  {ICCV} 2021, Montreal, QC, Canada, October 10-17, 2021}, pages 4563--4572.
  {IEEE}, 2021.

\bibitem{DBLP:conf/iclr/MengHSSWZE22}
Chenlin Meng, Yutong He, Yang Song, Jiaming Song, Jiajun Wu, Jun{-}Yan Zhu, and
  Stefano Ermon.
\newblock Sdedit: Guided image synthesis and editing with stochastic
  differential equations.
\newblock In {\em The Tenth International Conference on Learning
  Representations, {ICLR} 2022, Virtual Event, April 25-29, 2022}.
  OpenReview.net, 2022.

\bibitem{APF}
Michael~K. Pitt and Neil Shephard.
\newblock Filtering via simulation: auxiliary particle filters.
\newblock {\em J. Amer. Statist. Assoc.}, 94(446):590--599, 1999.

\bibitem{Filter_Finance}
B.~Ramaprasad.
\newblock {\em Stochastic filtering with applications in finance}.
\newblock 2010.

\bibitem{DBLP:journals/pami/SahariaHCSFN23}
Chitwan Saharia, Jonathan Ho, William Chan, Tim Salimans, David~J. Fleet, and
  Mohammad Norouzi.
\newblock Image super-resolution via iterative refinement.
\newblock {\em {IEEE} Trans. Pattern Anal. Mach. Intell.}, 45(4):4713--4726,
  2023.

\bibitem{DBLP:conf/iclr/SavinovCBEO22}
Nikolay Savinov, Junyoung Chung, Mikolaj Binkowski, Erich Elsen, and
  A{\"{a}}ron van~den Oord.
\newblock Step-unrolled denoising autoencoders for text generation.
\newblock In {\em The Tenth International Conference on Learning
  Representations, {ICLR} 2022, Virtual Event, April 25-29, 2022}.
  OpenReview.net, 2022.

\bibitem{Sny}
C.~Snyder, T.~Bengtsson, P.~Bickel, and J.~Anderson.
\newblock Obstacles to high-dimensional particle filtering.
\newblock {\em Mon. Wea. Rev.}, 136:4629--4640, 2008.

\bibitem{DBLP:conf/icml/Sohl-DicksteinW15}
Jascha Sohl{-}Dickstein, Eric~A. Weiss, Niru Maheswaranathan, and Surya
  Ganguli.
\newblock Deep unsupervised learning using nonequilibrium thermodynamics.
\newblock In Francis~R. Bach and David~M. Blei, editors, {\em Proceedings of
  the 32nd International Conference on Machine Learning, {ICML} 2015, Lille,
  France, 6-11 July 2015}, volume~37 of {\em {JMLR} Workshop and Conference
  Proceedings}, pages 2256--2265. JMLR.org, 2015.

\bibitem{NEURIPS2019_3001ef25}
Yang Song and Stefano Ermon.
\newblock Generative modeling by estimating gradients of the data distribution.
\newblock In H.~Wallach, H.~Larochelle, A.~Beygelzimer, F.~dAlch\'{e}-Buc, E.~Fox, and R.~Garnett, editors, {\em Advances in Neural
  Information Processing Systems}, volume~32. Curran Associates, Inc., 2019.

\bibitem{song2021scorebased}
Yang Song, Jascha Sohl-Dickstein, Diederik~P Kingma, Abhishek Kumar, Stefano
  Ermon, and Ben Poole.
\newblock Score-based generative modeling through stochastic differential
  equations.
\newblock In {\em International Conference on Learning Representations}, 2021.

\bibitem{10.1162/NECO_a_00142}
Pascal Vincent.
\newblock A connection between score matching and denoising autoencoders.
\newblock {\em Neural Comput.}, 23(7):1661–1674, jul 2011.

\bibitem{DBLP:conf/cvpr/WhangDTSDM22}
Jay Whang, Mauricio Delbracio, Hossein Talebi, Chitwan Saharia, Alexandros~G.
  Dimakis, and Peyman Milanfar.
\newblock Deblurring via stochastic refinement.
\newblock In {\em {IEEE/CVF} Conference on Computer Vision and Pattern
  Recognition, {CVPR} 2022, New Orleans, LA, USA, June 18-24, 2022}, pages
  16272--16282. {IEEE}, 2022.

\bibitem{DBLP:conf/icml/YuXMJPGZZW22}
Peiyu Yu, Sirui Xie, Xiaojian Ma, Baoxiong Jia, Bo~Pang, Ruiqi Gao, Yixin Zhu,
  Song{-}Chun Zhu, and Ying~Nian Wu.
\newblock Latent diffusion energy-based model for interpretable text modelling.
\newblock In Kamalika Chaudhuri, Stefanie Jegelka, Le~Song, Csaba
  Szepesv{\'{a}}ri, Gang Niu, and Sivan Sabato, editors, {\em International
  Conference on Machine Learning, {ICML} 2022, 17-23 July 2022, Baltimore,
  Maryland, {USA}}, volume 162 of {\em Proceedings of Machine Learning
  Research}, pages 25702--25720. {PMLR}, 2022.

\bibitem{zakai}
M.~Zakai.
\newblock On the optimal filtering of diffusion processes.
\newblock {\em Z. Wahrscheinlichkeitstheorie und Verw. Gebiete}, 11:230--243,
  1969.

\end{thebibliography}

\appendix

\section{Additional information on the implementation of diffusion models}\label{app:diff}
\subsection{The choice of the coefficients for the forward SDE}\label{app:fsde}
The task of the forward SDE in Eq.~\eqref{eq:forward} is to transform any given initial distribution $Q_0(Z_0)$ to the standard Gaussian distribution $\mathcal{N}(0, \mathbf{I}_d)$. It is shown in \cite{song2021scorebased,10.1162/NECO_a_00142,NEURIPS2020_4c5bcfec} that such task can be done by a linear SDE with properly chosen drift and diffusion coefficients. For example, we can define $b(\tau)$ and $\sigma(\tau)$ in Eq.~\eqref{eq:forward} by
\begin{equation}\label{eq:cof}
\begin{aligned}
b(\tau) = \frac{{\rm d} \log \alpha_\tau}{{\rm d} \tau} \;\;\; \text{ and }\;\;\; \sigma^2(\tau) = \frac{{\rm d} \beta_\tau^2}{{\rm d}\tau} - 2 \frac{{\rm d}\log \alpha_\tau}{{\rm d}\tau} \beta_\tau^2,
\end{aligned}
\end{equation}
where the two processes $\alpha_\tau$ and $\beta_\tau$ are defined by
\begin{equation}\label{eq:ab}
\alpha_\tau = 1-\tau, \;\; \beta_\tau = \tau \;\; \text{ for } \;\; \tau \in [0,1].
\end{equation}
The definitions in Eq.~\eqref{eq:cof} and Eq.~\eqref{eq:ab} can ensure that the conditional density function $Q_\tau(Z_\tau | Z_0)$ for any fixed $Z_0$ is the following Gaussian distribution:
\begin{equation}\label{eq:gauss}
Q_\tau(Z_\tau | Z_0) = \mathcal{N}(\alpha_\tau Z_0, \beta_\tau^2 \mathbf{I}_d).
\end{equation}
It is easy to see that the choice of $\alpha_\tau$ and $\beta_\tau$ can ensure that 
\[
Q_1(Z_1 | Z_0)  = \mathcal{N}(0, \mathbf{I}_d) \;\;\Longrightarrow\;\; Q_1(Z_1) = \int_{\mathbb{R}^d} Q_1(Z_1 | Z_0) Q_0(Z_0) dZ_0 = \mathcal{N}(0, \mathbf{I}_d),
\]
which is the property we need for the forward SDE.

\subsection{Discretization of the forward and reverse-time SDEs}\label{app:discrete}
Taking the reverse-time SDE in Eq.~\eqref{DM:RSDE} as an example, 
we use the Euler-Maruyama scheme to discretize the reverse-time SDE and transform any set of Gaussian samples $\{z_{1,j}\}_{j=1}^J$ of the final state $Z_1$ to a set of samples, denoted by $\{z_{0,j}\}_{j=1}^J$, approximately following the target distribution $Q_0(Z_0)$. Specifically, we first introduce a partition of the pseudo-temporal domain $\mathcal{T} = [0,1]$, i.e.,
\[
\mathcal{D}_K : = \{\tau_k \;| \;0 = \tau_0 < \tau_1  < \cdots < \tau_k < \tau_{k+1} < \cdots < \tau_K = 1\}
\] 
with uniform step-size $\Delta \tau = \f{1}{K}$. For each sample $z_{1, j}$, we obtain the approximate solution $z_{0,j}$ by recursively evaluating the following scheme
\begin{equation}\label{Euler:RSDE}
z_{\tau_k, j} = z_{\tau_{k+1}, j} - \big[ b(\tau_{k+1} ) z_{\tau_{k+1}, j}  - \sigma^2(\tau_{k+1}) S(z_{\tau_{k+1}, j}, \tau_{k+1}) \big] \Delta \tau + \sigma(\tau_{k+1}) \Delta W_{\tau_{k+1},j}, 
\end{equation}
for $k = K-1, K-2, \cdots, 1,0$, where $\Delta W_{\tau_{k+1},j}$ is a realization of the Brownian increment. The accuracy of $\{z_{0,j}\}_{j=1}^J$ is determined by the number of pseudo-time steps $K$.

\subsection{The loss function for training the diffusion model}\label{app:loss}
We provide details of the general loss function in Eq.~\eqref{optimization:score} for training the approximate score function. The full definition of the loss function in Eq.~\eqref{optimization:score} is 
\begin{equation}\label{eq:full_loss1}
Loss = \mathbb{E}_{\tau \sim \mathcal{U}[0,1],\, Z_0 \sim Q_0(Z_0),\, Z_\tau \sim Q_\tau(Z_\tau | Z_0)} \left[\lambda(\tau) \beta_\tau^2 
\| \nabla_z \log Q_\tau(Z_\tau) - \bar{S}(Z_\tau, \tau;\theta)\|_2^2\right],
\end{equation}
where $\mathcal{U}[0,1]$ is the uniform distribution, $Q_0(Z_0)$ is the target distribution, $Q_\tau(Z_\tau | Z_0)$ is the conditional distribution given in Eq.~\eqref{eq:gauss}, $\beta_\tau$ is defined in Eq.~\eqref{eq:ab}, and $\lambda(\tau)$ is a weighting function. This formulation is not practical because we do not know the exact score function $\nabla_z \log Q_\tau(Z_\tau)$. Thanks to the derivation in \cite{10.1162/NECO_a_00142}, the loss in Eq.~\eqref{eq:full_loss1} is equivalent to 
\begin{equation}\label{eq:full_loss2}
Loss = \mathbb{E}_{\tau \sim \mathcal{U}[0,1],\, Z_0 \sim Q_0(Z_0),\, Z_\tau \sim Q_\tau(Z_\tau | Z_0)} \left[\lambda(\tau) \beta_\tau^2 
\| \nabla_z \log Q_\tau(Z_\tau | Z_0) - \bar{S}(Z_\tau, \tau; \theta)\|_2^2\right] + const,
\end{equation}
where the exact score function is replaced by the gradient of the logarithm of the conditional distribution $Q_\tau(Z_\tau | Z_0)$. It makes the task much easier because we know $Q_\tau(Z_\tau | Z_0)$ is the Gaussian distribution $\mathcal{N}(\alpha_\tau Z_0, \beta_\tau^2 \mathbf{I}_d)$. Thus, we have 
\[
\nabla_z \log Q_\tau(Z_\tau | Z_0) = - \f{Z_\tau - \alpha_\tau Z_0}{\beta^2_\tau},
\]
such that the loss function in Eq.~\eqref{eq:full_loss2} is computable, i.e., 
\[
Loss = \mathbb{E}_{\tau \sim \mathcal{U}[0,1],\, Z_0 \sim Q_0(Z_0),\, Z_\tau \sim Q_\tau(Z_\tau | Z_0)} \left[\lambda(\tau)  
\left\| - \f{Z_\tau - \alpha_\tau Z_0}{\beta_\tau} - \beta_\tau\bar{S}(Z_\tau, \tau; \theta)\right\|_2^2\right] + const.
\]
Moreover, because $\f{Z_\tau - \alpha_\tau Z_0}{\beta_\tau}$ follows the standard Gaussian distribution, the final loss used to solve the optimization problem in Eq.~\eqref{optimization:score} is 
\[
Loss = \mathbb{E}_{\tau \sim \mathcal{U}[0,1],\, Z_0 \sim Q_0(Z_0),\, \zeta \sim \mathcal{N}(0, \mathbf{I}_d)} \left[\lambda(\tau)  
\left\| \zeta- \beta_\tau\bar{S}(Z_\tau, \tau; \theta)\right\|_2^2\right].
\]

\end{document}